\input amstex 
\documentstyle{amsppt}
\input bull-ppt

\font\tencyr=wncyr10
\def\cyr#1{\hbox{\tencyr#1}}

\keyedby{bull414/pxp}
\define\Spec{\operatorname{Spec}}
\define\Br{\operatorname{Br}}
\define\A{\operatorname{A}}
\define\PGL{\operatorname{PGL}}
\define\NS{\operatorname{NS}}
\define\cont{\operatorname{cont}}
\define\discrete{\operatorname{discrete}}
\define\loc{\operatorname{loc}}
\define\alg{\operatorname{alg}}
\define\gp{\operatorname{gp}}
\define\Inf{\operatorname{Inf}}
\define\Res{\operatorname{Res}}
\define\Hom{\operatorname{Hom}}
\define\Lim{\operatorname{Lim}}
\define\ord{\operatorname{ord}}
\define\Pic{\operatorname{Pic}}
\define\Tr{\operatorname{Tr}}
\define\N{\operatorname{N}}
\define\Gal{\operatorname{Gal}}
\define\Aut{\operatorname{Aut}}
\define\GL{\operatorname{GL}}
\define\SL{\operatorname{SL}}

\topmatter
\cvol{29}
\cvolyear{1993}
\cmonth{July}
\cyear{1993}
\cvolno{1}
\cpgs{14-50}
\title On the passage from local to global in number 
theory\endtitle
\author B. Mazur\endauthor
\shortauthor{B. Mazur}
\shorttitle{The passage from local to global in number 
theory}
\address Department of Mathematics, Harvard University, 1 
Oxford Street, 
Cambridge, Massachusetts 02140\endaddress
\ml mazur\@math.harvard.edu\endml
\date October 1, 1992 and, in revised form, February 17, 
1993\enddate
\endtopmatter

\document
\heading Introduction\endheading
Would a reader be able to predict the branch of 
mathematics that is the subject
of this article if its title had not included the phrase 
``in Number Theory''? 
The distinction ``local'' versus ``global'', with various 
connotations, has
found a home in almost every part of mathematics, local 
problems being often a
stepping-stone to the more difficult global problems.

To illustrate what the earthy geometric terms local and 
global signify in
Number Theory, let us consider any Diophantine equation, 
say, 
$$
X^n+Y^n+Z^n=0\qquad (n>2)\.
$$
Determining the rational (or equivalently, integral) 
solutions of such
equations can be difficult! In contrast, to determine its 
solutions in
integers modulo $m$ (at least, for any fixed modulus 
$m>0)$ is a finite
problem---a comparatively easy problem. Moreover, this 
easy problem, thanks to
the Chinese Remainder Theorem, can be replaced by the 
problem of finding its
solutions modulo the prime powers $p^\nu$ dividing $m$. 
Given an integral
solution, i.e., $a^n+b^n+c^n=0$, for $a,b,c\in\Bbb Z$, by 
viewing $(a,b,c)$ as
a triple of integers modulo $p^\nu$ for any prime number 
$p$ and any $\nu\ge
0$, we get, of course, a solution modulo $p^\nu$. 

Ever since (implicitly) the work of Kummer and 
(explicitly) the work of Hensel,
we know that it is useful to note that an integral 
solution to such a
Diophantine equation yields more: for each prime number 
$p$, and for positive
integers $\nu=1,2,3,\dotsc,$ if we choose to view the 
integral solution,
successively, as a solution modulo $p^\nu$, we then get, 
for every prime number
$p$, a system of solutions modulo $p^\nu$ for each 
$\nu=1,2,\dotsc,$ the system
being ``coherent'' in the sense that reduction of the 
$(\nu+1)$\<st solution
modulo $p^\nu$ yields the $\nu$\<th solution. Of course, 
such a system is
nothing more than a solution in the projective limit $\Bbb 
Z_p\coloneq\Lim(
\Bbb Z/p^\nu\Bbb Z)$, the ring of $p$-adic integers $(\Bbb 
Z_p$ being
considered as a topological ring, given the ``profinite 
topology'', i.e., the
natural topology that it inherits as a projective limit of 
the finite discrete
rings $\Bbb Z/p^\nu\Bbb Z)$. 

For more flexibility, one can work in $\Bbb Q_p$, the 
field of fractions of $
\Bbb Z_p$ (most economically obtained from $\Bbb Z_p$ by 
adjoining the inverse
of the single element $p)$. The field $\Bbb Q_p$ contains 
$\Bbb Q$ and can be
regarded as the topological (in fact, metrized) field 
obtained by completing
the field $\Bbb Q$ with respect to the ``$p$-adic'' 
metric\footnote{Two
rational numbers are very close in the $p$-adic metric if 
their difference,
expressed as a fraction in lowest terms, has the property 
that its numerator is
divisible by a high power of $p$. The standardly 
normalized $p$-adic metric 
$|\ |_p$ is characterized by the fact that it is 
multiplicative, that
$|p|_p=1/p$, and that if $a$ is an integer not congruent 
to $0\bmod p$, then
$|a|_p=1$.}. Evidently then, if a polynomial equation with 
$\Bbb
Q$-coefficients fails to have a solution over some 
$p$-adic field $\Bbb Q_p$,
it cannot have a solution over $\Bbb Q$.

The formal relationship between the field $\Bbb Q$ and the 
collection of
completions $\Bbb Q_p$ (for $p$ ranging through all prime 
numbers) is closely
analogous to the relationship between the field $K$ of 
global rational
functions on, say, a smooth algebraic curve $S$ over $\Bbb 
C$ in the affine
plane, and the collection of fields $K_s$ of formal 
Laurent series with finite
order poles, centered at points $s\in S$. The distinction 
of {\it global\/} v.
{\it local\/} is evident here: any global meromorphic 
function $f\in K$ gives
rise to a collection of local data, namely, the collection 
of Laurent
expansions of $f$ at each of the points $s\in S$. Our 
affine plane curve $S$ is
not compact, and therefore the collection of fields 
$\{K_s\}f_{s\in S}$ does
not reflect ``the whole local story''. To get ``the whole 
story'', we usually
augment our collection $\{K_s\}$ by including the 
completions of $K$
corresponding to the ``missing'' points, the points at 
infinity. It is only
when one takes account of all the completions $\{K_s\}_s$ 
for all points $s$ in
the smooth compactification of $S$ that certain global 
constraints on the
collections of local data hold, the most basic being that 
the number of zeros
of a global function (zeros counted with appropriate 
multiplicity) is equal to
the number of its poles; equivalently, if $\{f_s\}_s$ is a 
collection of
Laurent series all coming from the same global meromorphic 
function $f$, then
$$
\sum\ord_s(f_s)=0,\tag1
$$
where $\ord_s(f_s)$ is the order of zero (or, if negative, 
pole) of the Laurent
series $f_s$ at $s$.

In analogy, it has long been understood (cf. \cite{W}) 
that there is one
``missing field'' which should be considered along with 
the collection of
$p$-adic completions of $\Bbb Q$; namely, the field $\Bbb 
R$, of real numbers,
i.e., the $Ur$-completion of $\Bbb Q$ (with respect to the 
standard absolute
value metric, which we denote by $|\ |_\infty$, to bring 
it notationally into
line with the $p$-adic metrics). The analogue of (1) is 
the following (easily
proved) formula\footnote{Exponentiating both sides of this 
equation, one gets
the somewhat more familiar ``product formula''.}: If 
$r\in\Bbb Q$, then
$$
\sum_p\log|r|_p=0,\tag2
$$
provided the summation subscript $p$ ranges over all prime 
numbers and (``the
infinite prime'') $\infty$.

One weakness in the analogy between the collection 
$\{K_s\}_{s\in S}$ for a
compact Riemann surface $S$ and the collection $\{\Bbb 
Q_p$, for prime numbers
$p$, and $\Bbb R\}$ is that the fields $K_s$ for all 
points $s$ are very much
alike; they are even isomorphic. In contrast, no manner of 
squinting seems to
be able to make $\Bbb R$ the least bit mistakable for any 
of the $p$-adic
fields, nor are the $p$-adic fields $\Bbb Q_p$ isomorphic 
for distinct $p$. A
major theme in the development of Number Theory has been 
to try to bring $
\Bbb R$ somewhat more into line with the $p$-adic fields; 
a major mystery is
why $\Bbb R$ resists this attempt so strenuously. 

Here are two types of questions in Number Theory one might 
want to pursue by
passing from local to global:

\subheading{Type (A). Questions about rational points} 
Given a Diophantine
equation or a system of Diophantine equations with 
coefficients in $\Bbb Q$,
when does knowledge about its rational solutions over the 
collection of local
fields $\Bbb Q_p$ for all prime numbers $p$, and over 
$\Bbb R$, give us some
palpable information about its solutions over $\Bbb Q$?

Here is one of the simplest examples of such a question:

For which classes of systems of Diophantine equations with 
coefficients in $
\Bbb Q$ is it the case that the existence of solutions 
over $\Bbb Q_p$ for all
prime numbers $p$ and over $\Bbb R$ {\it guarantees\/} the 
existence of a
solution over $\Bbb Q$?

The question of existence of solutions of systems of 
Diophantine equations over
$\Bbb R$ or over $\Bbb Q_p$ is ``certifiably easy''; at 
least it is a decidable
question in the sense of mathematical logic. In fact, all 
first-order theories
over $\Bbb R$ (or over $\Bbb Q_p)$ are decidable, this 
being a result due to
Tarski (and Ax-Kochen), a concise treatment of which can 
be found in
\cite{Coh}.  Whether or not the question of existence of 
solutions of systems
of Diophantine equations over $\Bbb Q$ is a decidable 
question is still up in
the air\footnote{Matijasevic has established the 
algorithmic undecidability of
the question of whether any given Diophantine equation has 
an integral solution
\cite{Mat}.}, but, as the example which we trotted out at 
the beginning of this
article was meant to convince, decidable or not, even 
specific polynomials with
coefficients in $\Bbb Q$ can require quite elaborate 
theory before we can say
that all their $\Bbb Q$-rational points are 
known\footnote{Anyone who doubts
this should try to read published accounts of the proof of 
``Fermat's Last
Theorem'' for exponent $n=37$.}.

\subheading{Type (B). Passing from knowledge about local 
``structures'' to 
knowledge about global ``structures''} There is a wealth 
of literature on
various aspects of questions of type (A)\footnote{For 
example, for a recent
article which contains, among other things, an extensive 
review and
bibliography concerning current activities concerning the 
Hasse principle and
``weak approximation'' (of local solutions by global ones) 
in some classes of
varieties, see \cite{C-T--S-D}.}.  One cannot, really, 
effectively ``separate'' 
questions of type (A) from questions of type (B), but the 
central aim of this
article is to discuss an exciting development (finiteness 
of certain
Tate-Shafarevich groups\footnote{We will be explaining 
this.}), which is
conveniently expressed in the language of (B).

In 1986 at the MSRI, Larry Washington told Karl Rubin of 
some new ideas of
Thaine, suspecting that they might be helpful in Rubin's 
work.  These ideas, as
formulated by Thaine\footnote{Thaine says that at least 
partial inspiration for
these ideas came from a close reading of a paper of 
Kummer.}, had been brought
to bear on questions concerning the ideal class groups of 
cyclotomic fields
\cite{Th}. Rubin refashioned Thaine's ideas to serve as a 
tool in the study of
the arithmetic of complex multiplication elliptic curves. 
By their means, Rubin
was able, among other things, to prove that some 
Tate-Shafarevich groups of
elliptic curves over number fields were finite. This 
provided us with the first
examples (in the number field case\footnote{There had 
already been some
function field examples \cite{Mi1--3}.}) where the 
finiteness of some
Tate-Shafarevich groups was actually proved, despite the 
fact that (all)
Tate-Shafarevich groups had been conjectured, three 
decades ago, to be finite.

Soon thereafter, working independently, V. A. Kolyvagin 
devised his approach to
these finiteness questions (Kolyvagin's method of ``Euler 
Systems'' resembles
the method of Thaine-Rubin but is more flexible and is 
applicable to a much
broader 
range of problems).

Kolyvagin's ideas (and those of Rubin and Thaine) 
represent something of a
revolution\footnote{The impressive number of otherwise 
inaccessible problems
that these methods successfully treat is one measure of 
the revolutionary
aspect of this work; the bibliography at the end of these 
notes gives a small
sample of the activity involved here. For example, thanks 
to Kolyvagin and
Rubin we now have: (1) a significantly shorter proof of 
the ``Main Conjecture
over $\Bbb Q$'' (originally proved by Andrew Wiles and 
me); (2) the one-
and two-variable ``Main Conjecture'' over quadratic 
imaginary fields; and 
(3) (most of) the Birch Swinnerton-Dyer Conjecture for 
modular elliptic curves
whose $L$ function does not have a multiple zero at $s=1$, 
and in particular,
that (for these elliptic curves) the Tate-Shafarevich 
group is finite (see
Theorem 3).\endgraf
But even more telling is that, after all this, we now 
think of the subject in
quite a different way, and a range of new arithmetic 
problems seems nowadays to
be within reach. For example, there is Nekov\'a\v r's 
recent study of Selmer
groups attached to higher (even) weight to gain 
information about algebraic
cycles on Kuga-Sato varieties \cite{Ne1}; there is Flach's 
finiteness theorem
for certain Selmer groups attached to automorphic forms 
arising as the
symmetric square of classical modular forms which yields 
information about the
deformation space of the Galois representations attached 
to classical forms
\cite{F}; there is Darmon's attack on the ``refined 
conjectures'' of Birch
Swinnerton-Dyer type \cite{D1, D1} (see also \cite{BD}); 
there is Kato's recent
exciting announcement concerning the construction of 
Kolyvagin-type (``Euler'')
systems of classes in algebraic $K$-groups.}. 

My plan in Part I is to try to explain to a general 
audience some of the
importance of Finiteness of the Tate-Shafarevich Group. 
Anything discussed in
Part I that smacks of being too technical is tucked into 
footnotes or into Part
III.

Part II is an analytic continuation of Part I in which I 
try to give a bit of
the flavor of how one goes about showing such finiteness. 
It is written in more
technical language than Part I, but still no proofs are 
given here. I hope that
it may be useful to some readers as an accompaniment to 
some of the published
articles on the subject or as a companion to some of the 
other surveys; see,
for example, the recent \cite{Ne2}.

Part III provides proofs of some things that seemed 
natural to bring up in an
expository treatment of this sort, results that are, very 
likely, well known to
the experts but for which there is no reference in the 
literature.

\heading Contents\endheading

\noindent
Part I
\roster
\item "1." Local-to-global principles (an example of a 
question of type (A):
zeroes of quadratic forms)
\item "2." Local-to-global principles (an example of a 
question of type (B):
isomorphisms of projective varieties)
\item "3." First examples: projective space, quadrics, and 
curves of genus
zero
\item "4." Curves of genus $\ge 1$.  The ``local-to-global 
principle up to
finite obstruction''
\item "5." General projective varieties
\item "6." The Tate-Shafarevich group
\item "7." Implications of the Tate-Shafarevich Conjecture 
in the direction of
the general local-to-global principle (up to finite 
obstruction)
\item "8." Examples of finiteness of the Tate-Shafarevich 
group
\item "9." Returning to Selmer's curve
\endroster

\noindent
Part II
\roster
\item "10." The Tate-Shafarevich group and the Selmer group
\item "11." Local control of the 
Selmer group
\item "12." Class field theory, duality, and Kolyvagin 
test classes
\item "13." Rational points in extension fields yielding 
``test classes''
\item "14." Heegner points
\endroster

\noindent
Part III
\roster
\item "15." The cohomology of locally algebraic group 
schemes
\item "16." Discrete locally algebraic group schemes
\item "17." Finiteness theorems
\item "18." The mapping of the Tate-Shafarevich group of 
$A/K$ to $\scr S(A/K)$
\item "19." The local-to-global principle for quadrics
\endroster

\noindent
Bibliography

\heading Part I\endheading
\subheading{1. Local-to-global principles (an example of a 
question of type
(A)\footnote{Type (A) refers to the distinction alluded to 
in the discussion in
the introduction.}: zeros of quadratic forms)} The Theorem 
of Hasse-Minkowski
\cite{H} guarantees that 

\thm\nofrills{} a quadratic form with coefficients in 
$\Bbb Q$ admits a
nontrivial zero over $\Bbb Q$ if and only if it does so 
over $\Bbb Q_p$ for all
prime numbers $p$ and over $\Bbb R$.
\ethm

To get a closer view of the sort of information that this 
classical result
gives us, let us think about it in the following special 
case: Consider a
diagonal quadratic form in three variables
$$
F(X,Y,Z)=a\bfcdot X^2+b\bfcdot Y^2-c\bfcdot Z^2\tag 1
$$
with $a,b,c$ positive integers.

The fact that the coefficients in (1) are not all of the 
same sign guarantees
that our form has a nontrivial zero over $\Bbb R$. The 
Theorem of Chevalley
(which says that a polynomial in more variables than its 
degree, with integer
coefficients $\bmod p$ and no constant term has a 
nontrivial zero $\bmod p)$
gives us that (1) has a nontrivial zero modulo $p$ for any 
prime number $p$.
Moreover, for odd prime numbers $p$ not dividing $a\bfcdot 
b\bfcdot c$,
Hensel's Lemma then tells us that the nontrivial zero 
$\bmod p$ given to us by
Chevalley's theorem lifts to a nontrivial $p$-adic zero. 
Thus, in this case,
the Hasse-Minkowski Theorem assures us that, if 
$$
a\bfcdot X^2+b\bfcdot Y^2-c\bfcdot Z^2=0\tag 2
$$
has a nontrivial $p$-adic solution for the finite set of 
primes $p$ dividing 
$2\bfcdot a\bfcdot b\bfcdot c$, then it has a rational 
solution. 

There is no difficulty in stipulating congruence 
conditions on $a,b,c$ under
which (2) will have $p$-{\it adic solutions\/} for this 
finite set of prime
numbers $p$.

But even in a concrete instance of (2), where the 
Hasse-Minkowski Theorem
guarantees the existence of a rational solution, there 
remains much yet to
ponder about; for the proof of the Hasse-Minkowski Theorem 
does not produce any
specific nontrivial rational solution easily, despite the 
fact that the problem
of finding such solutions is perfectly effective (e.g., in 
Mordell's book
\cite{Mo} it is shown that if (2) has a nontrivial 
rational solution, it has
such a solution with
$$
|X|\le \sqrt{b\bfcdot c},\qquad |Y|\le \sqrt{c\bfcdot 
a},\qquad
|Z|\le\sqrt{a\bfcdot b}\.
$$
For other, more general bounds see \cite{Ca1, 6.8, Lemma 
8.1}).

As Victor Miller pointed out, it would be good to give 
explicit bounds for the
adelic version of this problem, i.e., given a finite set 
of prime numbers $p$
and specific $p$-adic points for these $p$, find rational 
points of prescribed
closeness to these $p$-adic points.

Moreover, it would be of interest to specify fast 
algorithms (e.g., are there
algorithms which are of {\it polynomial time\/} in 
$\max(\log a,\log b,\log
c)$?) for the determination of these rational points.

For recent results concerning the equidistribution 
properties of the sets $
\scr S_g(n)$ of integral solutions of $g(x,y,z)=n$ 
(for specified integers $n)$ where
$g(x,y,z)$ is a positive-definite ternary quadratic form 
(equidistribution of
the
sets $\scr S_g(n)$, that is, when they are scaled back to 
the unit
ellipsoid $g(x,y,z)=
1$ in $\Bbb R^3)$, see \cite{Du, GF, D--S-P}\footnote{The
case of three variables is particularly difficult; these 
equidistribution
results are based on delicate upper bounds for the 
absolute values of the
Fourier coefficients of modular forms of half-integral 
weight obtained by
Iwaniec \cite{I}.}. 

There are also distributional types of estimates coming 
from the circle method
concerning rational zeros of diagonal forms (in four or 
more variables). For
example, Theorem 4 of Siegel's 1941 paper \cite{Si} treats 
the special case of
a nondegenerate quadratic form in four variables with 
integral coefficients,
$$
F(X,Y,Z,T)=aX^2+bY^2+cZ^2+dT^2,
$$
whose determinant is the square of an integer. Explicitly, 
fix such a
quadratic form $F$, and for $S=(X,Y,Z,T)$ a four-tuple of 
integers, put
$$
\|S\|\coloneq|a|\bfcdot X^2+|b|\bfcdot Y^2+|c|\bfcdot Z^2+
|d|\bfcdot T^2\.
$$
Theorem 4 of loc. cit. gives a positive constant $\kappa$ 
such that 
$$
\sum_S\exp(-\|S\|/\lambda)=\kappa\bfcdot\lambda\log\lambda+o
(\lambda\log\lambda)\quad\text{(as }\lambda\to \infty),
$$
the summation being extended over all (nontrivial) 
integral zeros $S$ of $F$,
assuming that $F$ has nontrivial integral zeros 
(equivalently, assuming that
$F$ is indefinite and has a nontrivial zero in $\Bbb Q_p$ 
for all prime numbers
$p)$.

\subheading{2. Local-to-global principles (an example of a 
question of 
type (B)\footnote{Type (B) refers to the distinction 
alluded to in the
discussion in the introduction.}: isomorphisms of 
projective varieties)} One
says that two quadratic forms (in the same number of 
variables) over a field
$K$ are equivalent if one can pass from one of the 
quadratic forms to the other
(and from the other to the one) by suitable linear 
substitutions of variables,
where the matrix describing the linear change of variables 
has coefficients in
$K$. Questions of equivalence of quadratic forms over 
local fields are quite
manageable; for example, by ``Sylvester's Law'', quadratic 
forms over $\Bbb R$
are determined by three integers (their rank, nullity, and 
index).

Here is a variant of the Hasse-Minkowski theorem (cf. 
\cite{Has, Si}) expressed
as the answer to a question of type (B):

\thm\nofrills{} Two quadratic forms with coefficients in 
$\Bbb Q$ which are 
equivalent over $\Bbb Q_p$ for all prime numbers $p$ and 
over $\Bbb R$ are
equivalent over $\Bbb Q$.
\ethm

In this type (B) spirit, let us prepare for a more supple 
and more general
local-to-global question.

Let $V$ be a projective variety defined over any field 
$k$. We can think of $V$
as ``given by'' some finite collection of homogeneous 
forms in $N$ variables
(for some $N)$ with coefficients drawn from $k$. The set 
$V(L)$, of $L$-{\it
rational points\/} of $V$, for $L$ any field extension of 
$k$, is the common
``locus of zeros'' (in projective $N$-space with  
homogeneous coordinates in
$L)$ of the collection of homogeneous forms defining $V$. 
Assume that $V$ is
smooth\footnote{Why not? ``Smoothness'' is defined, as in 
multivariable
calculus, by the requirement that the appropriate jacobian 
determinant does not
vanish at $\overline k$-rational points of $V$.}. 

The notion of isomorphism between two (smooth) varieties 
$V$ and $W$ over $k$
(the technical term being {\it biregular 
isomorphism\/}\footnote{See Lecture 2
of \cite{Harr} for the notion of ``regular morphism'' over 
algebraically closed
fields and \cite{Hart} for the foundations for the more 
general concept of
``morphism of schemes''.} {\it defined over\/} $k)$ will 
be important to us.
But here, in contrast to the notion of equivalence of 
quadratic forms, one
allows changes of variables more fluid than the mere 
linear changes of
variables allowed in the definition of equivalence in the 
theory of quadratic
forms.\footnote{Nevertheless, in the special case of 
(smooth) quadrics, i.e.,
varieties given as the locus of zeros of a single 
(nondegenerate) quadratic
form, the quadratic forms are similar (i.e., equivalent up 
to multiplication by
a nonzero scalar) over $K$ if and only if the varieties 
are isomorphic over
$K$.  The point here is that (setting aside the case of 
quadric surfaces, for
the moment), if $V$ is a smooth quadric of dimension $1$ 
or $\ge 3$, then
$\Pic(V)=\Bbb Z$. Therefore, any automorphism of the 
variety $V$ must preserve
the isomorphism class of $L$, the (unique) generating 
ample line bundle, and
consequently must extend to projective linear 
automorphisms of the projective
spaces coming from the linear spaces of global sections of 
$L^{\otimes n}$ for
any $n$. Since the initial embedding of $V$ in projective 
space (via the
defining quadratic form) is the projective space of global 
sections of
$L^{\otimes 2}$ if $\dim V=1$, and of $L$ if $\dim V\ge 
3$, one easily sees
that similarity between such quadratic forms is the same 
as isomorphism between
their corresponding varieties. The same conclusion holds 
for smooth quadric
surfaces $V$ with a slightly different argument (since the 
Picard group of $V$
over $\overline K$ is $\Bbb Z\oplus \Bbb Z)$.} We admit 
the possibility that
$V$ and $W$ may be given as loci of zeros of systems of 
homogeneous forms in
different dimensional projective spaces and that the 
systems for $V$ and for
$W$ may comprise different numbers of homogeneous forms 
and of different
degrees.\footnote{When the field $k$ is $\Bbb C$, using 
the faithfulness of the
functor which passes from a smooth projective varieties 
$V$ over $\Bbb C$ to
its underlying complex analytic manifold $V(\Bbb C)$ 
\cite{Se1}, one sees that
to give a biregular morphism from $V$ to $W$ over $\Bbb C$ 
is equivalent to
giving a complex-analytic isomorphism from $V(\Bbb C)$ to 
$W(\Bbb C)$.}

Now let $V$ be a smooth variety over $\Bbb Q$. Embedding 
$\Bbb Q$ in its various
completions, we may forget that our $V$ is defined over 
$\Bbb Q$ and think of
$V$ as defined over $\Bbb Q_p$, where $p$ runs through all 
prime numbers, or
over $\Bbb R$, by simply regarding the defining collection 
of homogeneous forms
of $V$ as having their coefficients in $\Bbb Q_p$  or in 
$\Bbb R$. To indicate
that we have done this, we shall sometimes refer to the 
``original'' variety
$V$, defined over $\Bbb Q$ as $V/\Bbb Q$; and if we wish 
to think of it as
defined over $\Bbb Q_p$ or over $\Bbb R$, we refer to it 
as $V/\Bbb Q_p$ or $V/
\Bbb R$. 

\dfn{Question} To what extent do the ``local'' varieties
$$
\{V/\Bbb Q_p,\text{ for all }p,\text{ and }V/\Bbb R\}
$$
determine the ``global'' variety $V/\Bbb Q$ (up to 
isomorphism over $\Bbb Q)$?
\enddfn

To have vocabulary for this, let us say that a (smooth 
projective) variety $V'/
\Bbb Q$ is a {\it companion\/} to $V/\Bbb Q$ if their 
corresponding local
varieties are isomorphic,\footnote{{\it Companion\/} is a 
neologism; {\it
twisted\/} $\Bbb Q$-{\it form\/} is the term sometimes 
used, but I prefer a new
word here, since the term ``twisted form'' occurs in a 
number of different
contexts. Moreover, we will want to reserve the term 
``twist'' (below) for {\it
abelian varieties\/} which are twists of other abelian 
varieties.} i.e., if 
$V'/\Bbb Q_p$ is isomorphic to $V/\Bbb Q_p$ (as varieties 
over $\Bbb Q_p)$ for
all prime numbers $p$, and $V'/\Bbb R$ is isomorphic to 
$V/\Bbb R$ (as
varieties over $\Bbb R)$.

Denote by $\scr S(V)$ the set of isomorphism classes over 
$\Bbb Q$ of varieties
$V'/\Bbb Q$ which are companions to $V$.\footnote{So as 
not to get sidetracked
in technical things, we leave for Part III a discussion of 
the cohomological
interpretation of $\scr S(V)$.}

One can think of the cardinality of $\scr S(V)$ as roughly 
analogous to a {\it
class number\/}, i.e., a measure of the extent to which 
local data (in this
case, the isomorphism classes of $V/\Bbb Q_p$ for all $p$, 
and $V/\Bbb R)$
determine or fail to determine global data (the 
isomorphism class of $V/\Bbb
Q)$. One might say that the {\it local-to-global principle 
holds\/} for a class
of varieties $\scr V$ if $\scr S(V)$ consists of the 
single isomorphism class
$\{V\}$ for each member $V$ of $\scr V$.

\subheading{3. First examples: projective space, quadrics, 
and curves of genus
zero} 
\ch Projective space\endch 
The local-to-global principle holds for any variety $V/
\Bbb Q$ which is isomorphic (over $\Bbb C)$ to projective 
$N$-space
for some $N$. This statement,\footnote{Contained in the 
1944 Ph.D.  thesis of
Ch\^atelet, for a historical discussion of which see 
\cite{C-T, \S1}.} which is
one of the many manifestations of the {\it local-to-global 
principle for
elements of the Brauer group\/}\,\footnote{For the Brauer 
group and the
relationship between Brauer groups over local fields and 
global ones, see
\cite{Se3}.} {\it of\/} $\Bbb Q$, might be viewed as an 
arithmetic addendum to
the large body of literature (e.g., differential geometric 
and algebraic
geometric) establishing various ways in which projective 
space is {\it
rigid\/}. The analogous statements hold more generally 
with $\Bbb Q$ replaced
by any number field. 

\ch Quadrics\endch
The local-to-global principle holds for any smooth quadric 
variety $V$ over $
\Bbb Q$ (and, again, more generally over any number 
field). I am thankful to
Colliot-Th\'el\`ene for providing me with a sketch of a 
proof of this (see
Part III, \S19).

\ch Curves of genus $0$\endch
The local-to-global principle holds for any (smooth, 
projective) curve of genus
$0$ which is defined over $\Bbb Q$ (or over any number 
field). This follows
from {\it either\/} of the two previously cited examples, 
since:
\roster
\item "(i)" any smooth curve of genus $0$ over an 
algebraically closed field is
isomorphic to $\bold P^1$, and 
\item "(ii)" any smooth curve of genus zero over any field 
$k$ is canonically
isomorphic to a plane conic curve over $k$.\footnote{This 
isomorphism is given
as the projective embedding associated to the 
anticanonical line bundle on the
smooth curve of genus $0$. Serre suggested that I also 
point out the fact that
the category of quaternionic division algebras over any 
field $k$ of
characteristic $\ne 2$ is equivalent to the category of 
smooth curves of genus
$0$ over $k$, an equivalence in one direction being given 
as follows: If $D$ is
a quaternionic division algebra over $k$, let $V_D$ denote 
the plane conic over
$k$ given by the equations $\Tr(x)=0$ and $\N(x)=0$ in 
$D$, where $\Tr$ is the
trace and $\N$ is the quaternionic norm. {\it Exercise\/} 
(Serre). Find an
explicit equivalence of categories going in the other 
direction.}
\endroster

\subheading{4. Curves of genus $\ge 1$. The 
``local-to-global principle up to
finite obstruction''} That the local-to-global principle 
does not always hold for
curves of genus $1$ was seen early on, the first explicit 
examples being given
by Lind and Reichart in the early forties. In fact, this 
principle fails either
in its variants type (a) or type (B), for the failure of 
these two variants are
rather close in the context of curves of genus $1$: If a 
curve of genus $1$
possesses a $\Bbb Q_p$-rational point for all prime 
numbers $p$ and a real
point, and it does {\it not\/} have a $\Bbb Q$-rational 
point, 
then it has at least
one nonisomorphic companion. Conversely, if it has a 
nonisomorphic
companion, then at least one of its companions has no 
$\Bbb Q$-rational point
(yet has $\Bbb Q_p$-rational points for all $p$ and real 
points).

That the local-to-global principle fails even for smooth 
plane cubic
curves\footnote{These are all of genus $1$; smooth curves 
of degree $d$ in the
projective plane are of genus $g=(d-1)(d-2)/2$.} was shown 
by Selmer (1951).
Here is his most widely quoted example:

\thm\nofrills{}The curve $C\:3x^3+4y^3+5z^3=0$ has 
nonisomorphic 
companions. This
equation possesses nontrivial solutions over $\Bbb Q_p$ 
for all prime numbers
$p$ and over $\Bbb R$, but it possesses no nontrivial 
solutions over $\Bbb Q$.
\ethm

When I was preparing these notes, I initially intended to 
quote this Selmer
example and go on. But then it occurred to me that after 
the recent work which
is the subject of this article, we are actually in a 
position to explore the
question of companions to Selmer's curve (and others) with 
some precision.
Working out such an example seems to me to be a good way 
of demonstrating the
power of the recent results of Rubin and Kolyvagin. So, 
the following is an
explicit list of all the companions (worked out with the 
help of some e-mail
communications of Rubin):

\thm{Theorem 1} Selmer\RM's curve $C\:3x^3+4y^3+5z^3=0$ 
has, counting itself,
precisely five companions\,\RM:
$$
\align
3x^3+4y^3+5z^3=&0,\\
12x^3+y^3+5z^3=&0,\\
15x^3+4y^3+z^3=&0,\\
3x^3+20y^3+z^3=&0,\\
60x^3+y^3+z^3=&0\.
\endalign
$$
\ethm

\rem{Commentary} (1) This innocuous-sounding result lies 
quite deep. The fact
that all its companions are of degree 3 may be misleading; 
there is, a priori,
no guarantee that the companion curves to $V$ are 
isomorphic to smooth plane
curves over $\Bbb Q$, and, in fact, there is a priori no 
upper bound to their
degree in projective space. Establishing finiteness of 
$\scr S(C)$, something
we could not even begin to do before the work of Rubin and 
Kolyvagin, is
tantamount to establishing an upper bound to their degree.

(2) All five equations on this list have nontrivial 
rational solutions over $
\Bbb Q_p$ for all prime numbers $p$ and over $\Bbb R$. The 
first four equations
on the list possess no nontrivial rational solutions. The 
fifth equation
possesses a nontrivial rational solution $(0,1,-1)$, and 
this solution is
unique up to scalar multiplication (cf. \cite{Ca2, \S18}). 
If we take this
point as ``origin'' of the projective curve $E$ defined by 
the equation 
$$
60x^3+y^3=z^3=0,
$$
then $E$ is an elliptic curve over $\Bbb Q$ isomorphic to 
the jacobian of all
five curves on the list.

(3) In analogy with the discussion in \S2 of the 
Hasse-Minkowski Theorem,
giving the full list of companions as we have just done 
does not complete in a
satisfactory manner the discussion of {\it passage from 
local to global\/} for
this example: the demanding reader might ask how, if we 
are given some
projective smooth curve over $\Bbb Q$ together with the 
data that exhibits it
as a {\it companion\/} to Selmer's $C$, we may find a 
specific isomorphism over
$\Bbb Q$ between it and one of the members of our 
list.\footnote{I am thankful
for e-mail correspondence with Rubin in which he has 
proposed a direct and
succinct recipe for this. To put it in slightly different 
terms, suppose that
you are given a curve of $C$ of genus $1$ over $\Bbb Q$ as 
a projective variety
in $\bold P^N$ and you are merely told that $C$ has a 
rational point, but you
are also given a basis of the Mordell-Weil group of its 
jacobian $E$. Here is
how you can actually find a rational point on $C$. You 
first find any algebraic
point $P$ on $C$ of some degree $d$. Then find an 
algebraic point $\tau$ on $E$
such that $d\bfcdot \tau$ is the class of the divisor (of 
degree $0)$ on $C\:
\Sigma P^\sigma-d\bfcdot P$, where the summation is over 
all $\Bbb
Q$-conjugates of $P$. Let $\scr E$ denote a representative 
system for the group
of points on $E$ whose $d$\<th multiple lies in $E(\Bbb 
Q)$, modulo $E(\Bbb
Q)$. Then $\scr E$ is a finite set of points of $E$. 
Running through the
(finite) set of translates of $P$ by $e+\tau$ where $e$ is 
drawn from $\scr E$,
you are guaranteed that one of these is $\Bbb Q$-rational.}
\endrem

In view of this type of example, why not weaken somewhat 
the local-to-global
principle? Say that the {\it local-to-global principle 
holds, up to finite
obstruction\/} for a class of varieties $\scr V$, if each 
member $V$ of $\scr
V$ has only a finite number of nonisomorphic companions, 
i.e., if $\scr S(V)$
is a finite set for all $V$ of $\scr V$.

Despite the fact that we now, thanks to Rubin and 
Kolyvagin, can produce
examples of some curves of genus $1$ for which the 
local-to-global principle
holds up to finite obstruction, this question has not yet 
been resolved for all
curves of genus $1$. In contrast, consider:

\ch Curves of genus $\ge 2$\endch Here the local-to-global 
principle holds up
to finite obstruction; moreover, this is a rather easy 
result (see the
discussion in \S5 and Part III).

\ch Summary \endch The one outstanding unresolved case 
remaining for curves is
the case of genus $1$.

\subheading{5. General projective varieties}
\thm{Conjecture 1} The local-to-global principle holds, up 
to finite
obstruction for all \RM(smooth\RM) projective 
varieties\,\RM; that is, any
projective variety over $\Bbb Q$ has only a finite number 
of nonisomorphic
companions.
\ethm

\rem{Notes} As will be discussed in Part III, the set 
$\scr S(V)$, of
companions to $V$, depends only upon the nature of the 
$\Gal(\overline{\Bbb Q}/\Bbb Q)$-module of automorphisms 
of $V/\overline{\Bbb
Q}$, where $\overline{\Bbb Q}$ is an algebraic closure of 
$\Bbb Q$. For a
general projective variety, this $\Gal(\overline{\Bbb 
Q}/\Bbb Q)$-module
consists of the $\overline{\Bbb Q}$-rational points of a 
``locally algebraic
group'' $\underline{\Aut}(V)$.

In specific instances when the group $\underline{\Aut}(V)$ 
is easy to deal
with, we can prove Conjecture 1. For example, the proof of 
Conjecture 1 is
particularly easy if $\underline{\Aut}(V)$ is finite. This 
accounts for the
fact that we know the local-to-global principle up to 
finite obstruction for
{\it curves of genus\/}  $\ge 2$. For the same 
reason\footnote{Cf. \cite{KO}.},
we know the local-to-global principle up to finite 
obstruction for all
varieties $V/\Bbb Q$ of {\it general 
type\/}.\footnote{{\it General type\/}.
Let $r$ be a positive integer. Consider global tensors 
$\tau$ on a (smooth)
$n$-dimensional variety $V$ which for any system of local 
coordinates
$z_1,\dotsc,z_n$ can be expressed as $h(z)(dz_1\wedge 
dz_2\wedge\cdots\wedge
dz_n)^{\otimes r}$ where $h$ is some holomorphic function. 
If, for some choice
of $r$, there are ``enough'' linearly independent such 
tensors
$\tau_0,\dotsc,\tau_\nu$---enough in the sense that the 
rule which sends a
point $v\in V$ to the point in $\nu$-dimensional 
projective space with
homogeneous coordinates $[\tau_0(v),\dotsc,\tau_\nu(v)]$ 
is almost everywhere
defined and gives a birational imbedding of $V$ to 
$\nu$-dimensional projective
space---then one says that $V$ is of {\it general type\/}. 
The adjective
general signifies that, by some reckonings, many varieties 
have this property.
For example, for (smooth) curves, ``general type'' is 
equivalent to genus $>1$;
for surfaces, the varieties that are not of general type 
are either rational,
ruled, elliptic surfaces, abelian surfaces, $K3$ surfaces, 
or Enriques
surfaces, for smooth hypersurfaces $V$ of degree $d$ in 
projective $\nu$-space,
if $d\ge \nu+2$ then $V$ is of general type.} For almost 
the same reason the
local-to-global principle up to finite obstruction holds 
for all smooth {\it
hypersurfaces of dimension\/} $\ge 3$, or more generally, 
for all smooth {\it
complete intersections of dimension\/}  $\ge 3$ (see Part 
III, \S17, Corollary
4). 

A question, which arose in a conversation with N. Katz on 
these matters, is the
following: 
\ext
Is it within the range of the present ``state-of-the-art'' 
to show finiteness
of the part of $\scr S(V)$ represented by projective 
varieties whose numerical
invariants are ``limited'' in various ways? For example, 
fixing an integer $d$
define $\scr S(V)_d$ to be the subset of $\scr S(V)$ 
represented by projective
varieties which are companions to $V$ and which admit a 
projective imbedding of
degree $\le d$. Is $\scr S(V)_d$ finite for every $d$?
\endext
\endrem

\subheading{6. The Tate-Shafarevich group} In this 
section, we will be
considering abelian varieties over $\Bbb Q$; i.e., smooth, 
geometrically
connected, projective varieties over $\Bbb Q$ which are 
also algebraic groups,
the group law $A\times A\to A$ defined over $\Bbb Q$. Such 
algebraic groups are
necessarily abelian; since $A(\Bbb C)$, the complex points 
of $A$, form a
connected compact and commutative Lie group, $A(\Bbb C)$ 
is isomorphic, as a
topological group, to a product of circles.

If $A/\Bbb Q$ is an abelian variety, then its {\it 
Tate-Shafarevich group\/},
denoted 
$\cyr{X}(A/\Bbb Q)$, may be defined, first as a set, to be 
the set of
isomorphism classes over $\Bbb Q$ of companions of $A$ 
endowed with principal
homogeneous group action by $A$. More precisely, 
$\cyr{X}(A/\Bbb Q)$ is the set
of isomorphism classes over $\Bbb Q$ of pairs $(W,\alpha)$ 
where $W$ is a
projective variety defined over $\Bbb Q$ which is a 
companion of $A$ and $
\alpha\:A\times W\to W$ is a mapping of projective 
varieties, defined over $
\Bbb Q$, which is a principal homogeneous group action of 
$A$ on $W$. (Note.
To guarantee that the action $\alpha$ is a principal 
homogeneous action, one
need only check that it induces a principal homogeneous 
action of the Lie group
$A(\Bbb C)$ on the topological space $W(\Bbb C)$.)

The set $\cyr{X}(A/\Bbb Q)$ is given a (commutative) group 
structure via the
natural {\it Baer-sum construction\/} for principal 
homogeneous
spaces.\footnote{Baer-sum. If $W_1$ and $W_2$ are two 
principal homogeneous
spaces for $A$, let $A$ act on $W_1\times W_2$ by the 
antidiagonal action,
i.e., the action of $a$ in $A$ on $(w_1,w_2)$ is $(w_1+
a,w_2-a)$, where we have
written the action of $A$ on $W_1$ and on $W_2$ 
additively. Then the Baer-sum
of the principle homogeneous spaces $W_1$ and $W_2$ is the 
quotient variety of
$W_1\times W_2$ with respect to the antidiagonal action, 
this quotient variety
being viewed as principle homogeneous space for $A$ by an 
action induced from
the diagonal action of $A$ on $W_1\times W_2$.} Forgetting 
the principal
homogeneous action $(W,\alpha)\mapsto W$ gives a natural 
map (of pointed sets)
$$
\cyr{X}(A/\Bbb Q)\to \scr S(A/\Bbb Q)\.
$$
For an analysis of this mapping, see Part III, \S18, 
Theorem 5. In particular,
in \S17 we will see that there are a finite set $\scr A$ 
of abelian varieties
$A'/\Bbb Q$ and natural mappings 
$$
\cyr{X}(A'/\Bbb Q)\to\scr S(A/\Bbb Q)
$$
for each $A'\in\scr A$ such that the images of these 
mappings produce a
partition of the set $\scr S(A/\Bbb Q)$ indexed by $\scr 
A$ and the image of
$\cyr{X}(A'/\Bbb Q)$ in $\scr S(A/\Bbb Q)$ may be 
identified with the
orbit-space of $\cyr{X}(A'/\Bbb Q)$ under the natural 
action of the group of
$K$-automorphisms of $A'$. From this description (cf. 
explicitly, Part III,
\S18, Corollary 2 to Theorem 5) one sees that Conjecture 1 
restricted to the
class of abelian varieties over $\Bbb Q$ is equivalent to 

\thm{Conjecture 2 {\rm (Tate-Shafarevich)}\rm} For any 
abelian variety $A$ over
$\Bbb Q$, $\cyr{X}(A/\Bbb Q)$ is finite.
\ethm

The relation between $\cyr{X}(A/\Bbb Q)$ and $\scr 
S(A/\Bbb Q)$ has a
resonance, in the classical theory of Gauss and Lagrange, 
in the distinction
between the problems of classifying integral binary 
quadratic forms up to
strict equivalence and equivalence. The notion of 
equivalence of integral
binary quadratic forms is more intuitive than and is 
historically prior to the
notion of strict equivalence, but only with the latter 
notion does one get a
natural group structure on classes of forms.

What is known, in general, about the group structure of 
$\cyr{X}(A/\Bbb Q)$ is
that it is a torsion abelian group with the property that 
the kernel of the
homomorphism given by multiplication by any nonzero 
integer $n$ on it is
finite.

\subheading{7. Implications of the Tate-Shafarevich 
Conjecture in the direction
of the general local-to-global principle (up to finite 
obstruction)} In
preparing this article I have been tantalized by the urge 
to show that the
Tate-Shafarevich Conjecture (i.e., Conjecture 2 for 
abelian varieties over $
\Bbb Q)$ is equivalent to the local-to-global principle up 
to finite
obstruction (i.e., Conjecture 1) for all projective 
varieties over $\Bbb Q$.
This I have not done, but after some conversations and 
correspondence with
Yevsey Nisnevich and with Ofer Gabber, Theorem 2 has 
emerged. For its proof,
see Part III, \S17, Corollary 2 to Theorem 4. Before we 
state the result,
however, we must talk a bit about the group of connected 
components of
automorphism groups. 

Let $V/\Bbb Q$ be an arbitrary projective variety and 
$\underline{\Aut}(V)$ its
locally algebraic group of automorphisms. Let us denote by 
$\Gamma(V)$ the
group of connected components of the locally algebraic group
$\underline{\Aut}(V)/\overline{\Bbb Q}$, endowed with its 
natural continuous 
$\Gal(\overline{\Bbb Q}/\Bbb Q)$-action. Let $\Delta(V)$ 
denote the quotient of
$\Gal(\overline{\Bbb Q}/\Bbb Q)$ which acts faithfully on 
$\Gamma(V)$, and let
$\Gamma\rtimes \Delta$ denote the semidirect product 
constructed via this
action.

What is known about the structure of $\Gamma(V)$? As will 
be discussed in Part
III, there is a finitely generated abelian group $NS(V)$ 
with a continuous
$\Gal(\overline{\Bbb Q}/\Bbb Q)$-action on which the group 
$\Gamma(V)$ acts
with finite kernel (and in a manner compatible with its 
natural
$\Gal(\overline{\Bbb Q}/\Bbb Q)$-action). In particular, a 
quotient of
$\Gamma(V)$ by a finite group is isomorphic to a subgroup 
of $\GL(n,\Bbb Z)$
for some positive integer $n$, and therefore there is an 
upper bound to the
orders of finite subgroups of $\Gamma(V)$. It seems that 
little else is known
about the structure of $\Gamma(V)$! Is $\Gamma(V)$ 
finitely generated? Is it
finitely presented? Is it an arithmetic group? Conversely, 
what arithmetic
groups occur as $\Gamma(V)$\<'s?

One can show (Part III, \S17, Corollary 2):

\thm{Theorem 2} Let $V/\Bbb Q$ be a projective variety. 
Suppose that\RM:

{\rm (a)} the group $\Gamma(V)$ is finitely presented\RM;

{\rm (b)} there are only a finite number of distinct 
$\Gamma$-conjugacy classes
of finite subgroups in the semidirect product 
$\Gamma\rtimes\Delta$ 
\RM(i.e., the group $\Gamma(V)$ endowed with 
$\Gal(\overline{\Bbb Q}/\Bbb
Q)$-action is decent in the sense defined in Part {\rm 
III, \S16);} and 

{\rm (c)} the Tate-Shafarevich Conjecture 
holds\footnote{{\rm Or, somewhat more
specifically, that it holds for all twists of certain 
subabelian varieties of}
$\Pic^0(V)$ {\rm over} $\Bbb Q$.}.

Then $\scr S(V/\Bbb Q)$ is finite \RM(i.e., Conjecture 
{\rm 1} holds for $V)$.

\ethm

\rem{Some notes} If $\Gamma(V)$ is finitely generated, 
then $\Delta(V)$ is
finite. Ofer Gabber has sketched a proof that the 
conclusion of Theorem 2 holds
even if one weakens conditions (a) to the requirement that 
$\Gamma(V)$ be
finitely generated. Are there any projective varieties $V$ 
for which
$\Gamma(V)$ does not satisfy (a) and (b), i.e., for which 
$\Gamma(V)$ is not
decent?\footnote{Is it true that the automorphism groups 
of $K3$ surfaces, for
example, are decent? They are finitely generated. Are they 
arithmetic? I have
begun to pester experts about this. Relevant here are 
\cite{P-S--S, \S7,
Proposition and Theorem 1} which give somewhat explicit 
descriptions of the
automorphism groups of $K3$ surfaces: If $V$ is a $K3$ 
surface over
$\overline{\Bbb Q}$, $NS$ its N\'eron-Severi lattice, and 
$\Gamma$ its
automorphism group, then $\Gamma$ is commensurate with 
$O(NS)/W(NS)$ where
$O(NS)$ is the orthogonal group of the lattice $NS$, and 
$W(NS)$ is the
subgroup generated by reflections coming from elements in 
$NS$ of square $-2$.
For a finer study of which $\Gamma$\<'s occur, see Nikulin 
\cite{Nik1, 2}.}
\endrem

\subheading{8. Examples of finiteness of the 
Tate-Shafarevich group} Let $E/
\Bbb Q$ be a {\it modular\/} elliptic curve. That is, $E$ 
is an elliptic curve,
defined over $\Bbb Q$, whose underlying Riemann surface 
admits a nonconstant
holomorphic mapping from a compactification of the 
quotient of the upper half
plane under the action of a congruence subgroup of 
$\SL(2,\Bbb
Z)$.\footnote{This is, in fact, equivalent to the more 
usual formulation of
``modularity'', see \cite{Ma2}.} The conjecture of 
Taniyama-Weil asserts that
every elliptic curve over $\Bbb Q$ is modular. 

We must now introduce the $L$-function of $E/\Bbb Q$. For 
most of our purposes,
we may deal with a crude version of it, which we will call 
$L_S(E,s)$, where,
if we have fixed a model 
$$
y^2=x^3+ax+b\tag 3
$$
for $E$ over $\Bbb Q$, with $a,b\in\Bbb Z$, the subscript 
$S$ will be any
finite set of prime numbers which includes all the prime 
numbers dividing the
discriminant of this model. If $p$ is a prime number not 
in $S$, define the
{\it local factor\/} 
$$
L_p(E,s)\coloneq (1-a_pp^{-s}++p^{1-2s})^{-1},
$$
where $a_p$ is the integer defined by the formula $1+
p-a_p=$ the number of
rational points of the projective curve (over the prime 
field $\Bbb F_p)$
determined by reducing the above model modulo $p$.  Then 
$L_S(E,s)$ is defined
(initially) as the Dirichlet series
$$
\prod_{p\notin S}L_p(E,s),
$$
where $p$ runs through all prime numbers not in $S$. Using 
either elementary,
or not so elementary, estimates on the size of the $a_p$, 
one sees that the
Dirichlet series $L_S(E,s)$ converges to yield an analytic 
function of $s$ 
in a right half plane. Under the hypothesis that $E$ is 
modular, this analytic
function $L_S(E,s)$ extends to an entire function in the 
$s$-plane.

\rem{Remark} There is, in fact, a ``good way'' of 
extending the definition of
the {\it local factors\/} $L_p(E,s)$ for prime numbers $p$ 
in $S$.
The $L_p(E,S)$ are polynomials in $p^{-s}$ of degree $\le 
2$ (cf. \cite{Ta})
explicitly given in terms of a ``minimal model'' for the 
elliptic curve $E$
over $\Bbb Z$, and the ``good'' $L$ function $L(E,s)$ is 
defined to be the
product of these local factors $L_p(E,s)$ for all prime 
numbers $p$.

The methods of Kolyvagin, aided by results of Gross and 
Zagier, Waldspurger,
Bump, Friedberg, and Hoffstein, and R. Murty and K. Murty, 
yield the following
extraordinary result (which is an important piece of the 
Conjecture of Birch
and Swinnerton and Dyer).
\endrem

\thm{Theorem 3} If $E$ is a modular elliptic curve for 
which the Hasse-Weil $L$
function $L_S(E,s)$ either does not vanish at the point 
$s=1$ or has a simple
zero at $s=1$, then $\cyr{X}(E/\Bbb Q)$ is finite. 

In the former case, i.e., if $L_S(E,1)\ne 0$, then $E$ 
possesses only a finite
number of $\Bbb Q$-rational points. If $L_S(E,s)$ has a 
simple zero, then the
group $E (\Bbb Q)$ of $\Bbb Q$-rational points of $E$ is 
of rank $1$. 
\ethm

The order of zero of $L_S(E,s)$ at $s=1$ is independent of 
the choice of set of
primes $S$ and would be the same if we dealt with the 
``good'' $L$ function
$L(E,s)$. For refinements of the above theorem, see 
\cite{Coa, K1, K2, R1--8}.

\subheading{9. Returning to Selmer's curve} Briefly, this 
is what goes into the
proof of Theorem 1: the last of the five curves in the 
list (given in the
statement of Theorem 1) has a $\Bbb Q$-rational point and 
can be identified
with the jacobian, call it $E$, of all five curves (cf. 
\cite{Ca2, Chapters 18,
20} where it is shown that the cubic curve $ax^3+by^3+
cz^3=0$, for nonzero
rational numbers $a$, $b$, $c$ has, as jacobian, the 
elliptic curve
$abcx^3+y^3+z^3=0)$. This elliptic curve $E$ has complex 
multiplication, and a
computation gives that the ratio $L(E,1)/\Omega$ is equal 
to $9$, where
$L(E,s)$ is the (good) $L$-function of $E$ and $\Omega$ is 
the {\it real
period\/} of $E$. It follows then from the recent work of 
Rubin \cite{R5} that,
if $\cyr{X}=\cyr{X}(E/\Bbb Q)$ is the Tate-Shafarevich 
group of $E$ over $
\Bbb Q$, then $\cyr{X}$ is finite and of order equal to a 
power of 2 times a
power of 3. Checking prior computations of the 2- and 
3-primary components of
$\cyr{X}$ in tables of Stephens \cite{St}, one then has 
that $\cyr{X}$ is a
product of two cyclic groups of order 3. By Part III, 
\S18, Corollary 1 to
Theorem 5 the set of companions $\scr S(C)$ may be 
identified with the
quotient-set of $\cyr{X}$ under the involution $x\mapsto 
-x$ and, therefore,
has cardinality $=1+(9-1)/2=5$. One is left with the chore 
of finding
representatives for the 5 companions; this being done by 
the list presented in
the statement of the theorem, it remains only to check 
that all five curves on
our list possess rational points over $\Bbb Q_p$ for all 
prime numbers $p$, and
over $\Bbb R$, which they do, and that they are pairwise 
nonisomorphic over 
$\Bbb Q$, which they are\footnote{Here are some hints 
about this latter
assertion. Let $E[3]$ be the kernel of multiplication by 
$3$ in $E$, and note
that the natural mapping $H^1(\Bbb Q,E[3])\to H^1(\Bbb 
Q,E)$ is injective
(since $E(\Bbb Q)$ vanishes). Now if $C$ is any one of our 
five curves
$ax^3+by^3+cz^3$, let $\beta,\gamma\in\overline{\Bbb Q}$ 
be such that $\beta^3=b$, $\gamma^3=c$. The cohomology 
class corresponding to
$C$ in $H^1(\Bbb Q,E)$ is the image of the class in 
$H^1(\Bbb Q,E[3])$
represented by the 1-cocycle which sends an element 
$\sigma\in\Gal(\overline{
\Bbb Q}/\Bbb Q)$ to the point 
$(0,1,(\beta/\gamma)^{\sigma-1})$ in $E[3]$.
These 1-cocycles lie in distinct cohomology classes, for 
our five curves $C$.}.

\heading Part II\endheading
\subheading{10. The Tate-Shafarevich group and the Selmer 
group} Let 
$\overline{\Bbb Q}$ be the algebraic closure of $\Bbb Q$ 
in $\Bbb C$, and for
each prime number $p$ fix an embedding of $\overline{\Bbb 
Q}$ in $\overline{
\Bbb Q}_p$, an algebraic closure of $\Bbb Q_p$. Let $G$ 
denote the Galois group
$\Gal(\overline{\Bbb Q}/\Bbb Q)$, and for each $p$ let 
$G_p\coloneq 
\Gal(\overline{\Bbb Q}_p/\Bbb Q_p)$ and 
$G_\infty\coloneq\Gal(\Bbb C/\Bbb R)$, 
these being viewed as closed subgroups of $G$ (after our 
choices of embeddings).

Fix an elliptic curve $E/\Bbb Q$. If $K$ is a field 
extension of $\Bbb Q$, the
group of $K$-rational points of $E$ is denoted $E(K)$. 
Since $E$ is defined
over $\Bbb Q$, there is a natural (continuous) action of 
$G$ on $E(\overline{
\Bbb Q})$, of $G_p$ on $E(\overline{\Bbb Q}_p)$, and of 
$G_\infty$ on $E(\Bbb
C)$. The cohomological definition of $\cyr{X}(E/\Bbb Q)$ 
is as follows.
$$
\cyr{X}=\cyr{X}(E/\Bbb Q)\coloneq 
\bigcap\ker\{H^1(G,E(\overline{\Bbb Q}))\to
H^1(G_p,E(\overline{\Bbb Q}_p))\}
$$
where the intersection is taken over all ``finite'' primes 
$p$ and $p=\infty$.
It is useful to lighten the notation somewhat and to 
convene that
$$
\align
H^1(G,E)\coloneq&H^1(G,E(\overline{\Bbb Q})),\\
H^1(G_p,E)\coloneq&H^1(G_p,E(\overline{\Bbb Q}_p)),\\
H^1(G_\infty,E)\coloneq&H^1(G_\infty,E(\Bbb C))\.
\endalign
$$

If $A$ is an abelian group, and $n$ an integer, let $A[n]$ 
denote the kernel of
multiplication by $n$, so that we have the ``Kummer'' 
exact sequences of
continuous $G$-modules,
$$
0\to E(\overline{\Bbb Q})[n]\to E(\overline{\Bbb 
Q})\overset n\to\to
E(\overline{\Bbb Q})\to 0,
$$
and of continuous $G_p$-modules,
$$
0\to E(\overline{\Bbb Q}_p)[n]\to E(\overline{\Bbb 
Q}_p)\overset n\to\to
E(\overline{\Bbb Q}_p)\to 0,
$$
each yielding long exact sequences on cohomology, pieces 
of which, following
our convention, can be written as the horizontal exact 
sequences occurring in
the commutative diagram:
$$
\eightpoint
\CD
0 @>>> E(\Bbb Q)/n\bfcdot E(\Bbb Q) @>>> H^1(G,E[n])  
@>j>> H^1(G,E)[n]
@>>> 0\\
@.  @VVV   @VVV   @VVV\\
0 @>>>  E(\Bbb Q_p)/n\bfcdot E(\Bbb Q_p)  @>>> 
H^1(G_p,E[n]) @>>>
H^1(G_p,E)[n]  @>>> 0
\endCD
$$

Visibly, $\cyr{X}[n]$, the kernel of multiplication by $n$ 
in $\cyr{X}$, is the
intersection for all $p\le\infty$ of the kernels of the 
right-hand vertical
arrows in the above diagrams (for all $p)$.

\dfn{Definition} The $n$-{\it Selmer group\/} (of the 
elliptic curve $E/\Bbb
Q)$, denoted $S[n]$, is the subgroup of $H^1(G,E[n])$ 
defined as the full
inverse image of $\cyr{X}[n]\subset H^1(G,E)[n]$ under the 
map labelled $j$
above.
\enddfn

For any finite field extension $K$ of $\Bbb Q$ in 
$\overline{\Bbb Q}$, making
the base change from $\Bbb Q$ to $K$, we have the Tate 
Shafarevich group of
$E/K$ and the analogously defined $n$-Selmer group of 
$E/K$, which we will
denote $\cyr{X}_K$ and $S_K[n]$, respectively. Essentially 
by definition these
fit into an exact sequence
$$
0\to E(K)/n\bfcdot E(K)\to S_K[n]\to \cyr{X}_K[n]\to 0\.
$$

The importance of the Selmer group is that, as the above 
exact sequence
displays, it packages information concerning both the 
Mordell-Weil group and
the Tate-Shafarevich group; but the Selmer group is, in 
many instances,
curiously more tractable than either of them. This is 
rather like the
phenomenon met in the study of number fields, where, at 
times, the product of
{\it class number\/} times {\it regulator\/} is more 
readily computed than
either factor. In any event, to study elements of order 
$n$ in $\cyr{X}_K$, we
pass to the study of elements
$$
s\in S_K[n]\subset H^1(G_K,E[n]),
$$
where $G_K$ is the closed subgroup of $G$ which is the 
identity on $K$.

\subheading{11. Local control of the Selmer group} At this 
point we are ready
to discuss a step in Kolyvagin's program that may seem 
modest but is of crucial
importance in setting the stage for what is to come. The 
``abstract'' group
$E[n]$ is isomorphic to a product of two cyclic groups of 
order $n$ and
consequently its full automorphism group is isomorphic to 
the finite group 
$\GL_2(\bold Z/n\bold Z)$. The natural action of $G_K$ on 
$E[n]$ gives us, up
to equivalence, a representation $G_K\to \GL_2(\Bbb 
Z/n\Bbb Z)$, whose kernel
we will denote $G_M$, and, to be consistent, we use the 
letter $M$ to denote
the finite extension field of $K$ in $\overline{\Bbb Q}$ 
comprised of the
elements in $\overline{\Bbb Q}$ fixed by $G_M$. So, 
letting $\Delta_n$ denote
$\Gal(M/K)=G_K/G_M$, we have the natural inclusion 
$\Delta_n\subset\GL_2(\Bbb
Z/n\Bbb Z)$. 

Let us now say that the integer $n$ is good if the 
``restriction mapping'' on
cohomology, described below, is injective:
$$
H^1(G_K,E[n])\to 
H^1(G_M,E[n])^{\Delta_n}=\Hom_{\Delta_n}(G_M,E[n])\.
$$

The usefulness of this definition of goodness is that:

(a) invoking the known richness of action of the Galois 
groups of number fields
on $n$-torsion in elliptic curves, the hypothesis of 
goodness is not a strong
hypothesis on $n$ (and even when $n$ is not good, the 
kernel of the
``restriction mapping'', for $n$ replaced by powers of 
$n$, is finite, of order
bounded independent of the power).

(b) any element $h\in H^1(G_K,E[n])$ (and consequently, 
also, any element $s$
in the Selmer group) gives rise, by the ``restriction 
homomorphism'' above, to
a continuous $\Delta_n$-equivariant homomorphism; if $n$ 
is good, then
triviality of the $\Delta_n$-equivariant homomorphism 
$h\:G_M\to E[n]$ implies
triviality of the element $h\in H^1(G_K,E[n])$ giving rise 
to it. 

From now on we assume $n$ to be good, and we will freely 
identify elements of
$H^1(G_K,E[n])$ with the $\Delta_n$-equivariant 
homomorphisms $G_M\to E[n]$ to
which they give rise. This allows us to ``study'' elements 
$s\in S_K[n]$ by
local means, in the following sense.

Let $v$ run through all places of $K$, both nonarchimedean 
and archimedean, and
let $K_v$ denote the completion of $K$ at $v$, $\overline 
K_v$ an algebraic
closure containing $\overline{\Bbb Q}$, and 
$G_{K_v}\coloneq\Gal(\overline
K_v/K_v)$. Then we have the exact sequence
$$
0\to E(K_v)/n\bfcdot E(K_v)\to H^1(G_{K_v},E[n])\to 
H^1(G_{K_v},E)[n]\to 0,
$$
and, directly from the definition of the Selmer group, we 
have that the image
of an element $s\in S_K[n]$ under the natural homomorphism 
$H^1(G_K,E[n])\to
H^1(G_{K_v},E[n])$ maps to an element, call it $s_v$, in 
the subgroup
$$
E(K_v)/n\bfcdot E(K_v)\subset H^1(G_{K_v},E[n])\.
$$
Since $n$ is good, and since a homomorphism from $G_K$ is 
trivial if it is
trivial on all decomposition groups, we have something 
that could be called the
{\it principle of local control\/}:

\thm\nofrills{} If the {\rm ``}local\,{\rm ''} elements 
$s_v\in E(K_v)/n\bfcdot
E(K_v)$ vanish for all places $v$, then $s$ vanishes.
\ethm

\subheading{12. Class field theory, duality, and Kolyvagin 
test classes} The
$G_K$-module $E[n]$ admits a (nondegenerate) self-pairing, 
called the
``Weil-pairing''
$$
E[n]\times E[n]\to \mu_n
$$
where $\mu_n$ is the $G_K$-module of $n$\<th roots of 
unity. The Weil pairing
induces a pairing on local cohomology 
$$
\align
H^1(G_{K_v},E[n])\times H^1(G_{K_v},E[n])&\to 
H^2(G_{K_v},\mu_n)\subset\Bbb Q/
\Bbb Z\\
(x,y)&\mapsto \langle x,y\rangle_v
\endalign
$$
where the inclusion of $H^2(G_{K_v},\mu_n)$ in $\Bbb 
Q/\Bbb Z$ is given by
Local Class Field Theory. Important for us will be the 
fact\footnote{One might
add that the pairing $E(K_v)/n\bfcdot E(K_v)\times 
H^1(G_{K_v},E)[n]\to\Bbb Q/
\Bbb Z$ induced by the (Weil) self-pairing on 
$H^1(G_{K_v},E[n])$ is
compatible, up to sign, with the pairing coming from Tate 
Local Duality (for
abelian varieties). Cf. the related discussion \cite{Mi4, 
Chapter I, \S3,
pp.~54--55}.} that if both $x$ and $y$ map to zero under 
the homomorphism to
$H^1(G_{K_v},E)$, then $\langle x,y\rangle_v=0$.

The Weil pairing also induces a pairing on global 
cohomology 
$$
\align
H^1(G_K,E[n])\times H^1(G_K,E[n])&\to 
H^2(G_K,\mu_n)\subset\bigoplus_v
H^2(G_{K_v},\mu_n)\subset\bigoplus_v\Bbb Q/\Bbb Z\\
(x,y)&\mapsto\bigoplus\langle x_v,y_v\rangle_v
\endalign
$$
where $x_v$ and $y_v$ are the images of $x$ and $y$ in the 
local
cohomology
$H^1(G_{K_v},E[n])$, the direct sum on the right being 
taken over places $v$ of
$K$. The inclusion on the right is that given by Global 
Class Field Theory,
which embeds $H^2(G_K,\mu_n)$ into the kernel of the sum 
mapping $\bigoplus_v
\Bbb Q/\Bbb Z\to\Bbb Q/\Bbb Z$.

\ch Kolyvagin's basic strategy\endch 
With our data, $E$, $K$, $n$, fixed and understood, let 
$w$ be a place of $K$.
By a {\it Kolyvagin test class for\/} $w$ let us mean a 
cohomology class $c\in
H^1(G_K,E[n])$, which has the property that $c$ goes to 
$0$ under the
composition of natural mappings $H^1(G_K,E[n])\to 
H^1(G_{K_v},E[n])\to
H^1(G_{K_v},E)$ for all places $v\ne w$ of $K$ and does 
not go to $0$ for $v=w$.

Let us first see that any Kolyvagin test class imposes a 
local condition on
every element of the Selmer group. For, if $s\in S_K[n]$ 
and $c$ is a Kolyvagin
test class, we have that $\langle c_v,s_v\rangle_v=0$ for 
$v\ne w$, since for
places different from $w$, both $c_v$ and $s_v$ map to 
zero in
$H^1(G_{K_v},E)$.

On the other hand, by the consequence of Global Class 
Field Theory cited above,
we have
$$
\sum \langle c_v,s_v\rangle_v=0,
$$
the summation being taken over all places $v$ of $K$.

Putting these together, we get a {\it local condition\/} 
at $w$, which is
satisfied for any element $s$ in the Selmer group 
$S_K[n]$; namely, 
$$
\langle c_w,s_w\rangle_w=0\.\tag 1
$$

Kolyvagin's strategy, roughly put, is to produce a large 
and systematic
collection of {\it test classes\/} $c$ so as to obtain 
enough conditions of the
type (1) on the local components of elements $s$ to 
provide the tightest
bounds on the size of the Selmer group.

For example, in one series of applications, $n$ is a prime 
number $>2$, the
field $K$ is a quadratic imaginary field, and, letting the 
superscript ${}^\pm$
denote $\pm$ eigenspace for the nontrivial automorphism of 
$K$ (we are in a
situation where each of the two eigenspaces 
$\{E(K_w)/n\bfcdot E(K_w)\}^\pm$
is cyclic of order $n)$. In this setup, the mere existence 
of a {\it Kolyvagin
test class\/} for $w$ in the $\pm$ eigenspace guarantees 
(via condition (1))
that, for any Selmer class $s\in S_K[n]^\pm$ (same sign), 
we have $s_w=0$.

Naturally, one might (and in fact one does) wish, at 
times, to loosen the
requirement on a test class $c$ that goes to $0$ in 
$H^1(G_{K_v},E)$ for all
but one place $v$. Fixing a finite set $S$ of places, one 
can make similar good
use of classes $c$ that go to $0$ in $H^1(G_{K_v},E)$ for 
all $v$ outside the
finite set $S$. 

\subheading{13. Rational points in extension fields 
yielding ``test classes''}
The basic mechanism that links rational points to 
Kolyvagin test classes is
encapsulated by the diagram displayed below. Here, $L/K$ 
is a finite Galois
extension with Galois group $\scr G$, and we assume that 
$L/K$ is unramified
outside the single place $w$ of $K$ and that the 
``restriction mapping'',
labelled $\Res$ below, is an isomorphism:
$$
\eightpoint
\CD
@.  @.   @. 0\\
@.  @.   @.  @VVV\\
@.  @.   @.  H^1(\scr G,E)[n]\\
@.  @.   @.  @VV\Inf V\\
0  @>>>  E(K)/n\bfcdot E(K)  @>>> H^1(G_K,E[n])  @>>>  
H^1(G_K,E)[n]\\
@.  @VVV   {\ssize{\cong}}@VV\Res V   @VVV\\
0  @>>> \{E(L)/n\bfcdot  E(L)\}^{\scr G} @>>>  
\{H^1(G_L,E[n])\}^{\scr G}
@.\subset H^1(G_L, E[n])\,^{\scr G}
\endCD
$$

Now let $P\in E(L)$ be a rational point satisfying the 
following two
conditions---the first being the ``serious'' condition, 
and the second a minor
one:

(a) For all $\sigma\in\scr G$, $P^\sigma-P$ is divisible 
by $n$ in $E(L)$.

(b) For each place $\nu$ of $L$ which does not lie over 
the place $w$ of $K$,
the specialization of $P$ to the group of connected 
components of the N\'eron
fiber of $E$ at $\nu$ is of order relatively prime to $n$.

For N\'eron models, see Artin's exposition (in \cite{A2}) 
of N\'eron's original
paper \cite{N\'e}; see also Grothendieck's Exp. IX in 
\cite{Grot3} and
\cite{BLR}.

Let us call rational points $P\in E(L)$ satisfying (a), 
(b) above {\it
Kolyvagin rational points\/}.

By (a), $P$ gives rise to a class $[P]$ in 
$\{E(L)/n\bfcdot E(L)\}^{\scr G}$.

\dfn{Definition} The class $c=c_P\in H^1(G_K,E[n])$ is the 
{\it unique homology
class\/} such that $\Res(c)$ is equal to the image of 
$[P]$ in
$\{H^1(G_L,E[n])\}^{\scr G}$.
\enddfn

There is a unique homology class $d\in H^1(\scr G,E)[n]$ 
such that $\Inf(d)$ is
the image of $c$ in $H^1(G_K,E)[n]$. Using (b) one can 
show that the
restriction of the class $d$ to decomposition groups 
attached to every place
$v\ne w$ of $K$ vanishes, i.e., that $c_P$ is a {\it 
Kolyvagin test class\/}
for $w$ provided its image in $H^1(G_{K_w},E)$ does not 
vanish. Clearly a
necessary requirement for this nonvanishing to occur is 
that $[P]$ be a ``new''
class in $\{E(L)/n\bfcdot E(L)\}^{\scr G}$, i.e., a class 
not in the image of
$E(K)/n\bfcdot E(K)$. 

\ch Summary\endch
We have started with an elliptic curve $E/\Bbb Q$. We have 
fixed a ``good''
nonzero integer $n$. We have somewhat silently passed to a 
finite extension
field $K/\Bbb Q$. Let us call this $K$ the {\it first 
auxiliary field
extension\/}: as already hinted, it will be chosen to be a 
suitable quadratic
imaginary field extension of $\Bbb Q$. For any place $w$ 
of $K$, we now wish to
find a suitable Galois extension field $L/K$ unramified 
outside $w$ (call $L/K$
the {\it second auxiliary field extension\/}) and also a 
``Kolyvagin'' rational
point $P$ in $E(L)$ which yields a ``new'' class  $[P]$ in 
$\{E(L)/n\bfcdot
E(L)\}^{\scr G}$.

The construction of a suitable ``first'' auxiliary field 
extension $K/\Bbb Q$
and then an ample and systematic supply of ``second'' 
auxiliary field
extensions $L/K$ together with rational points $P$ is the 
punch line of
Kolyvagin's program. For then we would have a significant 
supply of Kolyvagin
test classes $c_P$ which would impose a significant number 
of conditions on the
local components of elements of the $n$-Selmer group. By 
the ``principle of
local control'', one then would have control of the 
$n$-Selmer group itself.

\subheading{14. Heegner points} For any quadratic 
imaginary field extension $K/
\Bbb Q$ (an appropriate one of these will be our choice of 
``first auxiliary
field extension'') and for any positive integer $f$, let 
$K_f$ denote the {\it
ring class field of conductor\/} $f$ {\it over\/} $K$ (for 
suitable $f$,
$K_f/K$ will be our choice of ``second auxiliary field 
extension'' $L/K)$.

We have a tower of field extensions, with Galois Groups as 
marked
 below, and where the full Galois group $\Gal(K_f/\Bbb Q)$ 
is equal to a
``generalized dihedral group'', i.e., an extension of a 
cyclic group of order 2
by the abelian group $\scr G_f\coloneq\Gal(K_f/K)$, with 
action of the
nontrivial element in the cyclic group on $\scr G_f$ given 
by multiplication by
$-1$. 
$$
\matrix
\format\l\\
K_f\\
\ |\ \}(\scr O_K/f\bfcdot\scr O_K)^*/(\Bbb Z/f\Bbb Z)^*\\
K_1\\
\ |\ \}\text{ the ideal class group of } K\\
K\\
\ |\ \}\text{ cyclic of order two}\\
\Bbb Q
\endmatrix
$$

Here, $\scr O_{K}$ denotes the ring of integers of $K$. If 
$\scr O_{K,f}$
denotes the {\it order\/} $\Bbb Z+f\bfcdot\scr 
O_K\subset\scr O_K$, we have that
$\Gal(K_f/K)$ is isomorphic to the group of isomorphism 
classes of locally free
$\scr O_{K,f\,}$-modules of rank one.

As previously mentioned, the Taniyama-Weil conjecture 
asserts that any elliptic
curve $E/\Bbb Q$ can be realized as a quotient curve (over 
$\Bbb Q)$ of the
modular curve $X_0(N)$ for a suitable positive integer 
$N$. Noncuspidal points
$x$ of $X_0(N)$ may be ``identified'' with pairs $(\scr 
E,C)$ up to
isomorphism, where $\scr E$ is an elliptic curve and $C$ 
is a cyclic subgroup
of $\scr E$ of order $N$. The point $x$ is rational over a 
given field
extension of $\Bbb Q$ if and only if the pair $(\scr E,C)$ 
is definable over
that number field. Suppose then that we have a nonconstant 
mapping, defined
over $\Bbb Q$, $\varphi\:X_0(N)\to E$.

We now choose a quadratic imaginary field $K\subset \Bbb 
C$ such that all prime
factors of $N$ split in $\scr O_K$. Under this hypothesis, 
there is an ideal $
\scr N\subset\scr O_K$ such that $\scr O_K/\scr N$ is 
cyclic of order $N$; we
choose such an ideal. 

For any positive integer $f$ relatively prime to $N$, the 
ideal $\scr
N_f\coloneq\scr N\cap\scr O_{K,f}$ is an invertible $\scr 
O_{K,f}$-module such
that $\scr O_{K,f}/\scr N_f$ is cyclic of order $N$. Let 
$\scr E_f$ be the
elliptic curve $\Bbb C/\scr O_{K,f}$, and let $C_f\subset 
\scr E_f$ be the
cyclic subgroup of order $N$ given by $\scr N_f^{-1}/\scr 
O_{K,f}$. Let
$x_{K,f}\in X_0(N)$ be the point corresponding to the pair 
$(\scr E_f,C_f)$,
and put $y_{K,f}=\varphi(x_{K,f})$.

Then $y_{K,f}$ is a $K_{f\,}$-rational point of $E$. We 
refer to it as the {\it
Heegner point of\/} $E$ (for the choices: $K\subset \Bbb 
C$, $\scr N$, $f$, and
the parametrization $\varphi)$. The collection of Heegner 
points for a fixed
$K$ and all positive integers $f$ which are relatively 
prime to $N$ satisfy the
axioms of what Kolyvagin calls a ``Euler system'', these 
axioms being the
essence of what is needed to produce a supply of 
``Kolyvagin points''.

To avoid getting at all into the technicalities of the 
passage from Heegner
points to Kolyvagin points, and yet to give some small, 
but honest, hint of the
flavor of it, consider the following:

We restrict attention to the case where $n$ is a prime 
number. Now choose $f$
to be a prime number which is {\it inert\/} (and 
unramified) in $K$, which does
not divide $N\bfcdot n$, and such that the Frobenius 
element at $(f\,)$ in
$G_{\Bbb Q}$ acts as a nonscalar involution on the group 
of $n$-torsion points
of the elliptic curve $E$. Then $\Gal(K_f/K_1)$ is cyclic 
of order $f+1$. Fix a
generator $\gamma\in\Gal(K_f/K_1)\subset \Gal(K_f/K)$, and 
fix a representative
system $S\subset\Gal(K_f/K)$ modulo 
$\Gal(K_f/K_1)$ so that every element of
$\Gal(K_f/K)$ is in $\gamma^i\bfcdot S$ for a unique 
integer $i$ in the range
$0\le i\le f$. Put
$$
P_{K,f}\coloneq \sum i\bfcdot\gamma^i\bfcdot s\bfcdot 
y_{K,f}\in 
E(K_f),
$$
the summation being taken over all $s\in S$ and $i$ in the 
above range.

Then (see \cite{Gros} for a very neatly written proof of 
this) taking our
``second auxiliary choice'' $L/K$ to be $K_f/K$, the point 
$P_{K,f}$ is a
``Kolyvagin point''; i.e., it satisfies (a), (b) above, 
and, moreover, the
class $[P_{K,f}]$ in $\{E(L)/nE(L)\}^{\scr G}$ is 
independent of the choice of
representative system $S$; it is independent up to scaling 
by $(\Bbb Z/n\Bbb
Z)^*$ of the choice of generator $\gamma$.

\heading Part III\footnote{{\rm I am very thankful to Ofer 
Gabber and Yevsey
Nisnevich for the crucial help they gave to me while I was 
writing
this.}}\endheading
\subheading{15.  The cohomology of locally algebraic group 
schemes} In Part I,
for ease of notation, the base field was taken to be $\Bbb 
Q$, but here it is
more natural to allow arbitrary number fields $K$ (of 
finite degree over $
\Bbb Q)$. For $v$ a place of $K$, we let $K_v$ denote the 
completion of $K$
with respect to $v$. So, if $V$ is a projective variety 
over $K$, we have the
analogous definition of {\it companion\/} variety to $V$, 
i.e., a projective
variety $V'$ over $K$ which is isomorphic to $V$ over 
every completion $K_v$.
We also have $\scr S(V/K)$, the set of isomorphism classes 
of companions (over
$K)$.

Let $X$ denote the spectrum of the ring of integers in 
$K$, and we reserve the
letter $S$ for finite sets of places of $K$, containing 
all infinite places.
Let $X-S$ be the open subscheme in $X$ which is the 
complement of the closed
points corresponding to the nonarchimedean places of $S$.

We shall be studying group schemes over $X-S$ which are 
(separated and) locally
of finite type but not necessarily of finite type. To 
signal this, we shall
refer to such group schemes as {\it locally algebraic 
group schemes\/}
(extending the terminology in \cite{BS} where locally 
algebraic groups are over
fields rather than over more general base schemes). We 
shall refer to a group
scheme which is of finite type as an {\it algebraic group 
scheme\/}.

Let $G$, then, be a smooth locally algebraic group scheme 
over $X-S$. That is,
$G$ is a group scheme (separated, and), locally of finite 
type over $X-S$, for
which $G^0$, the connected component containing the 
identity section, is a
smooth connected open subgroup scheme of $G$, of finite 
type over $X-S$. If
$Z\to X-S$ is any morphism of schemes, the pullback of our 
group scheme $G$ to
the base $Z$ is denoted $G/Z$. Whenever we write 
$H^1(Z,G)$ we can mean,
equivalently, \'etale or flat (fppf) cohomology of $G/Z$ 
over $Z$ (compare
\cite{Grot}). Given a $G$-valued 1-cocycle $c$ for the 
\'etale topology on $Z$,
we can use the conjugation action of $G$ on $G$ to 
``twist'' the locally
algebraic group scheme $G/Z$ by $c$ (an ``inner twist''; 
cf. \cite{Gi}) and the
isomorphism class of the resulting locally algebraic group 
scheme over $Z$ is
dependent only upon the cohomology class $\gamma$ in 
$H^1(Z,G)$ of our cocycle
$c$; we shall denote this ``inner twist'' by $G_\gamma/Z$.

Since our groups are, in the main, noncommutative, we are 
referring here to
noncommutative one-dimensional cohomology, given the 
structure only of pointed
set.  Giving the concepts  ``$\prod$\<'' and ``$\ker$\<'' 
the evident meanings in
the context of pointed sets, we let 
$$
\cyr{X}(G)\coloneq \ker\left\{H^1(K,G)\to \prod_v 
H^1(K_v,G)\right\},
$$
the product being taken over all places $v$ of $K$, and 
$$
\cyr{X}_S(G)\coloneq\ker\left\{H^1(X-S,G)\to \prod_{v\in 
S}H^1(K_v,G)\right\}\.
$$

\ex{Examples} (1) by a ``locally constant'' group scheme 
$G/S$ we mean a group
scheme which is locally constant for the \'etale topology, 
i.e., such that
there is a finite \'etale cover $S'\to S$ for which the 
pullback $G/S'$ is a
constant group scheme. Let $G/K$ be a {\it locally 
constant group\/} over the
number field $K$, and let $L/K$ be the (minimal) finite 
splitting field for the
action of $\Gal(\overline K/K)$ on $G(\overline K)$. Then 
$\cyr{X}(G)$ is
isomorphic to 
$$
\ker\left\{H^1(L/K,G)\to\prod_v H^1(K_v,G)\right\}\.
$$

The reason for this is the following. First, if the action 
of $\Gal(\overline
K/K)$ on $G(\overline K)$ is trivial, then $\cyr{X}(G)$ is 
trivial, because in
that case $H^1(K,G)$ is identified with conjugacy classes 
of continuous
homomorphisms from $\Gal(\overline K/K)$ to $G(\overline 
K)$, and if a
homomorphism restricts to the trivial homomorphism from 
$\Gal(\overline
K_v/K_v)$ to $G(\overline K_v)$, for all places $v$, then 
it is trivial. In the
general situation, then, an element of $\cyr{X}(G)\in 
H^1(K,G)$ must map to the
trivial element in $H^1(L,G)$ and therefore come from 
$H^1(L/K,G)$, as stated.
\endex

\thm{Corollary} Let $G$ be a locally constant group whose 
Galois action is
split by a finite cyclic extension of $K$. Then 
$\cyr{X}(G)$ is trivial.
\ethm

\demo{Proof} Let $L/K$ be the finite cyclic extension 
which splits the Galois
action of $G$. By the Cebotarev Theorem, there is a place 
$w$ of $L$ lying over
a place $v$ of $K$ (which is unramified, and) for which
$\Gal(L_w/K_v)=\Gal(L/K)$. It then follows that for this 
$v$, the mapping
$H^1(L/K,G)\to H^1(K_v,G)$ is injective.\qed\enddemo

(2) Let $G$ be a locally algebraic group over $K$ which is 
an
extension\footnote{We keep to the convention that a group 
$G$ is an extension
of $A$ by $B$ if $B$ is the normal subgroup of $G$ and $A$ 
is the quotient
group $G/B$.} of a locally constant group which is an 
arithmetic group (whose
continuous $\Gal(\overline K/K)$-action comes from an 
action of an ambient
linear algebraic group over $\Bbb Q)$ by a linear 
algebraic group (such groups
are said to be of ``type ALA'' in \cite{BS}). Then 
$\cyr{X}(G)$ is finite and,
moreover, the mapping
$$
H^1(K,G)\to\prod_vH^1(K_v,G)
$$
is proper, i.e., has finite fibers \cite{BS, Theorem 7.1}.

If $G$ is a connected, simply connected, semisimple 
algebraic group over $K$,
then this mapping is now known to be an isomorphism (if 
one either includes or
neglects to include the factors on the right corresponding 
to nonarchimedean
$v$, since $H^1(K_v,G)=0$ for finite $v$, and connected, 
simply connected,
semisimple algebraic groups $G)$. This fact, going under 
the heading {\it The
Hasse Principle\/} for connected, simply connected, 
semisimple algebraic
groups, has been the fruit of a long development (cf. 
\cite{Se2, Kne} for a
discussion of the work on this problem for the classical 
semi-simple groups;
Harder's paper \cite{Hard} for a proof of the above 
assertion for all
connected, simply connected, semisimple algebraic groups 
$G/K$ which do not
have a factor of type $E_8$; Platonov's survey \cite{Pl} 
for a neat discussion
of where the problem stood circa 1982; and Chernousov's 
resolution of what
remained to be done with respect to $E_8$ \cite{Ch}, 
thereby finishing the
problem).

\subheading{16. Discrete locally algebraic group schemes} 
Fix $\overline K$ an
algebraic closure of $K$, and consider the category of 
{\it discrete\/}
$\Gal(\overline K/K)$-{\it groups\/}, i.e., ``abstract 
groups'' (viewed as
discrete topological groups) endowed with continuous 
$\Gal(\overline
K/K)$-action. This category is equivalent to the category 
of discrete locally
algebraic groups over $K$, the equivalence 
$$
\discrete\,\loc\,\alg\,\gp\ \scr G/K 
\rightarrowtail\discrete
\Gal(\overline K/K)\text{-}\gp\ \Gamma
$$
being given by taking $\Gamma$ to be the group $\scr 
G(\overline K)$ with its
natural  $\Gal(\overline K/K)$ Galois action. The discrete 
locally algebraic
group $\scr G/K$ is locally constant (in the sense 
discussed in the previous
section, i.e., locally constant for the \'etale topology) 
if and only if the
action of $\Gal(\overline K/K)$ on $\Gamma$ factors 
through a finite quotient
group (this being automatic if $\Gamma$ is finitely 
generated).

Let $\Gamma$ be a discrete $\Gal(\overline K/K)$-group. 
Denote by
$\Delta=\Gal(L/K)$ the quotient of $\Gal(\overline K/K)$ 
which acts faithfully
on $\Gamma$; i.e., $L\subset \overline K$ is the splitting 
field of the action
on $\Gamma$. Let $\Gamma\rtimes\Delta$ denote the 
semidirect product
constructed via the action of $\Delta$ on $\Gamma$. 

If $\Gamma$ is finitely generated, then $\Delta$ is 
finite. In this case, if $S$
denotes a (finite) set of places including all archimedean 
places and all
places ramified in the finite extension $L/K$, and if 
$\scr G/K$ is the
discrete, locally algebraic group over $K$ associated to 
$\Gamma$, then $\scr
G/K$ is the generic fiber of a (unique) locally algebraic 
group scheme $\scr
G/X-S$ which is locally constant for the \'etale topology 
over $X-S$.

\dfn{Definition} A discrete $\Gal(\overline K/K)$-group 
$\Gamma$ (and its
associated discrete, locally algebraic group $\scr G/K)$ 
is called {\it
decent\/} if 
\roster
\item "(a)" the group $\Gamma$ is finitely presented;
\item "(b)" the group $\Gamma\rtimes\Delta$ (notation as 
above) has only a
finite number of $\Gamma$-conjugacy classes of finite 
subgroups.
\endroster
\enddfn

We will also call a discrete locally algebraic group 
scheme over $X-S$ 
{\it decent\/} if its generic fiber is so.

A sub-\<$\Gal(\overline K/K)$-group of finite index in a 
decent $\Gal(\overline
K/K)$-group is again decent. An inner twist of a decent 
group is decent. An
arithmetic group $\Gamma$ endowed with trivial Galois 
action is decent.

\thm{Lemma (a)} If $\scr G=\scr G/X-S$ is a decent 
\RM(discrete\RM) locally
constant group scheme, then $H^1(X-S,\scr G)$ is finite 
and \RM(therefore\RM)
so is $\cyr{X}_S(\scr G)$.
\ethm

\demo{Proof} Let $\Gamma$ be the associated discrete 
$\Gal(\overline
K/K)$-group. Letting $Y$ be the spectrum of the ring of 
integers of $L$, the
splitting field of the natural Galois action on $\Gamma$, 
and $T$ the inverse
image of $S$ in $Y$, we have that $Y-T$ is a ``Galois'' 
extension of $X-S$,
with group $\Delta=\Gal(L/K)$. We may suppose that $Y-T$ 
is a finite \'etale
Galois extension of $X-S$.

The set $H^1(Y-T,\scr G)$ is finite, for we may identify 
it with the
orbit-space of $\Hom_{\cont}(\pi_1(Y-T),\Gamma)$ under the 
action of
$\Gamma$-conjugation, where $\Hom_{\cont}$ means 
``continuous homomorphism''.
Since the image of $\pi_1(Y-T)$ in $\Gamma$ under any 
continuous homomorphism
is a finite subgroup of $\Gamma$, by ``decency'' of 
$\Gamma$, we have that
there are at most a finite number of $\Gamma$-conjugacy 
classes of finite
subgroups of $\Gamma$, and our finiteness assertion follows.

Consequently there is a finite \'etale Galois extension 
$Y'-T'\to X-S$
admitting $Y-T\to X-S$ as an intermediate extension, such 
that every class in
$H^1(Y-T,\scr G)$ is trivialized in $Y'-T'$. Denote by 
$\Delta'$ the Galois
group of $Y'-T'$ over $X-S$ and by 
$\psi\:\Delta'\to\Delta$ the natural
surjection. Then $H^1(X-S,\scr 
G)=H^1(\Delta',\Gamma)=\{\Gamma$-conjugacy
classes of homomorphisms 
$\varphi\:\Delta'\to\Gamma\rtimes\Delta$ which lift
$\psi\:\Delta'\to\Delta\}$. 
Since, by the decency hypothesis, the image of $\varphi$ 
has only a finite
number of possibilities up to $\Gamma$-conjugacy, so does 
$\varphi$.\qed\enddemo

\subheading{17. Finiteness theorems} In this section we 
let $G$ be a smooth
locally algebraic group scheme over $X-S$ and $G^0\subset 
G$ the open subgroup
scheme which is the ``connected component'' of $G$ 
containing the identity
section. The quotient $G/G^0$ viewed as sheaf for the 
\'etale topology over
$X-S$ is representable, in general, only as an {\it 
algebraic space group\/}
over $X-S$ (i.e., group object in the category of 
algebraic spaces over $X-S$;
cf. \cite{A1, Knu}); even if it is representable as a 
group scheme, the group
scheme does not have to be separated. Denote by $\scr 
G/X-S$ the algebraic
space group $G/G^0$ over the base $X-S$. As usual, we let 
$\scr G/K$ denote its
generic fiber viewed as (discrete, locally algebraic) a 
group over $K$.

From now on we shall assume that two further properties 
hold for $G$:

(A) The algebraic space group $\scr G/X-S=G/G^0$ is 
representable as a {\it
locally constant group\/} (for the \'etale topology) over 
$X-S$ and is decent
(cf. \S16).

(B) The smooth algebraic group scheme $G^0/X-S$ is an 
extension of an abelian
scheme $A/X-S$ by a connected linear affine group scheme.

The second property above is assumed only to tidy up the 
mental picture we have
of our locally algebraic group scheme: it is not really 
necessary for what we
will do and can always be achieved by suitable 
augmentation of the finite set
$S$. Under our assumptions, then, our locally algebraic 
group scheme $G$ is a
successive extension of the following three basic 
``building blocks'':
\roster
\item "(a)" a discrete, decent, locally constant group 
scheme $\scr G/X-S$,
\item "(b)" an abelian scheme $A/X-S$,
\item "(c)" a connected (smooth) linear affine group 
scheme $B/X-S$.
\endroster

We refer below to the subquotient $A$ of $G$ as the {\it 
abelian scheme part\/}
of $G$. Let $A^0\subset A$ denote the ``connected 
component'' of the N\'eron
model of $A$, i.e., the open subgroup scheme all fibers of 
which are connected.

\thm{Lemma (b)} If $A/\Bbb Q$ is an abelian variety and 
$S$ contains all the
primes of bad reduction of $A/\Bbb Q$ \RM(so that $A$ 
extends to an abelian
scheme over $X-S)$, then $\cyr{X}_S(A/X-S)$ is isomorphic 
to the image of
$H^1(X,A^0)$ in $H^1(X-S,A)$. The group $\cyr{X}_S(A/X-S)$ 
is finite if
$\cyr{X}(A)$ is finite.
\ethm

\demo{Proof} One has a commutative diagram, in which the 
horizontal rows are
exact:
$$
\eightpoint
\CD
@. H^1(X,A^0) @>>> H^1(X-S,A^0)  @>>> \prod\limits_{v\in 
S}H^1(K_v,A^0)\\
@.  @. @VV=V     @VV=V\\
0  @>>>  \cyr{X}_S(A/X-S) @>>> H^1(X-S,A)  @>>> 
\prod\limits_{v\in S}H^1(K_v,A)
\endCD
$$
the middle arrow being an equality because all primes of 
bad reduction for $A$
are assumed to be in $S$. Exactness of the top horizontal 
row can be seen, for
example, from the discussion in the Appendix of \cite{Ma1, 
p.~263}; exactness 
of the bottom horizontal row is simply by the definition 
of $\cyr{X}_S(A/X-S)$.
It follows that there is a natural surjection
$$
\alpha\:H^1(X,A^0)\twoheadrightarrow\cyr{X}_S(A/X-S)%
\subset H^1(X-S,A)\.
$$

Let $\Sigma$ denote the image of $H^1(X,A^0)$ in 
$H^1(X,A)$. By the appendix of
\cite{Ma1} $\Sigma$ contains $\cyr{X}(A)$ as a subgroup of 
finite index (in
fact, this index is a power of two). By the above 
discussion, there is a
natural surjection of $\Sigma$ onto $\cyr{X}_S(A/X-S)$.

This proves the lemma.\qed\enddemo

\thm{Lemma (c)} If $B/X-S$ is a connected \RM(smooth\RM) 
linear affine group
scheme, then $H^1(X-S,B)$ is finite.
\ethm

\demo{Proof} This is \cite{Nis1, Theorem 3.7}; compare 
also the discussion
around \cite{Nis2, Theorem 3.10}.\qed\enddemo

To deal with groups which are extensions of these building 
blocks, consider an
exact sequence of locally algebraic smooth group-schemes 
over a scheme $Y$,
$$
1\to U\to G\overset p\to\to W\to 1\.\tag 1
$$
If $\gamma\in H^1(Y,G)$ is a one-dimensional 
(noncommutative) cohomology class,
let $U_\gamma$ be the group scheme $U$ twisted by the 
class $\gamma$, where $G$
acts on $U$ by conjugation. Let 
$$
\pi\:H^1(Y,G)\to H^1(Y,W)
$$
be the mapping induced from $p$ on one-dimensional 
(noncommutative) cohomology.
If $\omega\in H^1(Y,W)$, let $H^1(Y,G)_\omega$ denote the 
fiber of $\pi$ above
$\omega$. If $\gamma\in H^1(Y,G)$ and 
$\omega=\pi(\gamma)$, by the ``long exact
sequence'' for noncommutative cohomology \cite{Gi, Chapter 
III,
3.3.5}\footnote{For the Galois-cohomological version, 
compare \cite{Se2,
Chapter VII, Appendix, Proposition 2}.}, there is a 
natural surjection of sets 
$$
\sigma\:H^1(Y,U_\gamma)\to H^1(Y,G)_\omega
$$
(where the image of the trivial element in 
$H^1(Y,U_\gamma)$ is the element $
\gamma$ in $H^1(Y,G)_\omega)$.

Now suppose that $Y=X-S$. We have the diagram 
$$
\CD
1 @>>>  \cyr{X}_S(G)  @>>> H^1(X-S,G)  @>>> 
\prod\limits_{v\in S} H^1(K_v,G)\\
@.  @VVV  @V\pi VV  @VV\pi V\\
1 @>>> \cyr{X}_S(W)  @>>> H^1(X-S,W)  @>>> 
\prod\limits_{v\in S}H^1(K_v,W)
\endCD
$$

\thm{Extension Lemma} Given the exact sequence $(1)$, 
suppose that

{\rm (i)} $\cyr{X}_S(W)$ is finite\RM;

{\rm (ii)} The $\cyr{X}_S(U_\gamma)$ are all finite 
\RM(for all twisting classes
$\gamma\in H^1(Y,G));$\newline 
and either

{\rm (iii1)} The group scheme $U$ is commutative and the 
maps
$H^1(K_v,U_\gamma)\to H^1(K_v,G)_\omega$ are all proper 
(for all $v\in S$, and
all twisting classes $\gamma\in 
H^1(X-S,G)$, where
$\omega=\pi(\gamma))$,\newline
or

{\rm (iii2)} The $H^1(X-S,U_\gamma)$ are all finite 
\RM(for all twisting
classes $\gamma\in 
H^1(X-S,G))$.\newline
Then $\cyr{X}_S(G)$ is finite.
\ethm

\demo{Proof} If $\omega$ is an element in the image of 
$\cyr{X}_S(G)$ in 
$\cyr{X}_S(W)$ and if $\gamma\in\cyr{X}_S(G)\subset 
H^1(X-S,G)$ is such that $
\omega=\pi(\gamma)$, consider the diagram
$$
\CD
1 @>>> \cyr{X}_S(U_\gamma)  @>>> H^1(X-S,U_\gamma)  @>>> 
\prod\limits_{v\in
S}H^1(K_v,U_\gamma)\\
@.  @VVV  @V\sigma VV  @VVV\\
1 @>>> \cyr{X}_S(G)_\omega  @>>> H^1(X-S,G)_\omega @>>> 
\prod\limits_{v\in
S}H^1(K_v,G)
\endCD
$$
where $\cyr{X}_S(G)_\omega$ is defined so as to make the 
lower sequence exact.

Since, by (i), $\cyr{X}_S(W)$ is finite, our lemma will 
follow if we show
$\cyr{X}_S(G)_\omega$ to be finite for each element 
$\omega$ in the image of
$\cyr{X}_S(G)$ in $\cyr{X}_S(W)$. Fix, then, such a 
$\omega\in\cyr{X}_S(W)$,
and an element $\gamma\in\cyr{X}_S(G)$ such that 
$\omega=\pi(\gamma)$.

From (ii) and (iii1) we see that the ``diagonal mapping'' 
of the right-hand
rectangle
$$
H^1(X-S,U_\gamma)\to\prod_{v\in S}H^1(K_v,G)
$$
has finite kernel. By surjectivity of $\sigma$ we see that 
this kernel maps
onto $\cyr{X}_S(G)_\omega$ proving the lemma in the case 
when (iii1) holds.

Now suppose (iii2), giving us finiteness of 
$H^1(X-S,U_\gamma)$. But
since
$\sigma\:H^1(X-S,U_\gamma)\to H^1(X-S,G)_\omega$ is 
surjective, finiteness of
$\cyr{X}_S(G)_\omega$ is immediate.\qed\enddemo

\thm{Theorem 4} Let $G$ be a locally algebraic 
\RM(smooth\RM) group scheme over
$X-S$ satisfying our two properties {\rm (A)} and {\rm 
(B)}. Suppose, further,
that the Tate-Shafarevich Conjecture holds for abelian 
varieties over $K$. Then
$\cyr{X}_S(G)$ is finite.
\ethm

\rem{Remark} If $S\subset S'$ is an inclusion of finite 
sets of places of $K$,
and if $G$ is a locally algebraic (smooth) group scheme 
over $X-S$ satisfying 
properties (A) and (B), then the restriction of $G$ to 
$X-S'$ also satisfies
these properties, so we might amplify the conclusion of 
our theorem to say that
$\cyr{X}_{S'}$ also satisfies these properties, so we 
might amplify the
conclusion of our theorem to say that $\cyr{X}_{S'}(G)$ is 
finite for $S'$ any
finite set of places of $K$ containing $S$.
\endrem

\demo{Proof} We do this in two steps.

First we shall assume that our locally algebraic group 
scheme $G$ fits into an
exact sequence of locally algebraic group schemes over 
$X-S$, 
$$
1\to U\to G\overset p\to\to W\to 1\tag1
$$
where $W$ is a discrete, decent, locally constant group 
scheme and $U$ is an
abelian scheme. Let us show that, in this situation, 
$\cyr{X}_S(G)$ is finite.

For this, we use the Extension Lemma, noting that property 
(i) holds by Lemma
(a), property (ii) holds by Lemma (b), and our assumption 
of finiteness of the
Tate-Shafarevich groups (of abelian varieties which are 
``twists'' of $A)$. We
now show that property (iii1) holds. We use the long exact 
sequence of
noncommutative cohomology groups and are significantly 
helped here by the fact
that $U$ is an abelian scheme and, in particular, 
commutative. Properness of
the mapping $\varphi\:H^1(K_v,U_\gamma)\to 
H^1(K_v,G)_\omega$ (for all twisting
classes $\gamma)$ is equivalent to showing that the 
kernels are finite. Let
$L_w/K_v$ be a finite local field extension over which the 
discrete locally
constant group $W$ is constant, i.e., such that the action 
of $\Gal(\overline
L_w/L_w)$ on $W$ is trivial. We have a commutative square
$$
\CD
H^1(K_v,U_\gamma) @>\varphi>> H^1(K_v,G)_\omega\\
@VVV    @VVV\\
H^1(L_w,U_\gamma)  @>\varphi'>> H^1(L_w,G)_\omega
\endCD
$$
and a coboundary mapping
$$
\delta\:H^0(L_w,W)\to H^1(L_w,U_\gamma)
$$
which is a homomorphism of groups and whose image is the 
kernel of $\varphi'$.
To see that $\ker\varphi$ is finite, it suffices to show 
that the kernel of the
diagonal mapping $\psi\:H^1(K_v,U_\gamma)\to 
H^1(L_w,G)_\omega$ in the square
above is finite. Since the kernel of the left-vertical map
$H^1(K_v,U_\gamma)\to H^1(L_w,U_\gamma)$ is finite, 
finiteness of $\psi$ will
follow if we show that the image of the homomorphism 
$\delta$ is finite. But
since $W$ is a constant decent group over $L_w$, the 
domain of $\delta$,
$H^0(L_w,W)$, is finitely generated. The range of $\delta$ 
is an abelian
torsion group, and therefore the image of $\delta$ is 
finite. This concludes
step 1.

Now we suppose that our locally algebraic groups scheme 
$G$ is an extension
(1) where $W$ is any locally algebraic group scheme with 
$\cyr{X}_S(W)$ is
finite, and where $U$ is a connected linear affine group 
scheme. In this case,
by assumption we have (i) and by Lemma (c) we have (iii2) 
(which already
implies (ii)) so that the Extension Lemma applies again to 
yield our
theorem.\qed\enddemo

\thm{Corollary 1} Let $G$ be a locally algebraic group 
over $K$ whose group of
connected components is decent. Suppose, further, that the 
Tate-Shafarevich
Conjecture holds. Then $\cyr{X}(G)$ is finite.
\ethm

\demo{Proof} We begin with a lemma on prolongations of 
locally algebraic 
groups over $K$ to (smooth) locally algebraic group 
schemes over $X-S$ 
(some $S)$.\enddemo

\thm{Lemma} Let $G$ be a locally algebraic group over $K$ 
whose group of
connected components $\scr G/K$ is decent. Then $G$ admits 
an extension to a
smooth locally algebraic group scheme, $G/X-S$, satisfying 
{\rm (A)} and 
{\rm (B)} over $X-S$, for some finite set of places $S$.
\ethm

\demo{Proof} Let $G^0/K$ denote the connected component of 
the identity in
$G=G/K$, so that $\scr G/K=G/G^0$ is decent, by our 
hypothesis. We must, for a
suitable finite set of places $S$, construct a 
``prolongation of $G/K$'', i.e.,
a smooth locally algebraic group scheme $G/X-S$ with the 
property that its
generic fiber is $G/K$ and its quotient by $G^0/X-S$, the 
connected component
in $G/K$ containing the identity, is locally constant 
(with generic fiber
isomorphic to $\scr G)$. We then can achieve (A) and (B) 
by further
augmentation of $S$.

Let $\Gamma=\scr G(\overline K)$ be the associated 
discrete $\Gal(\overline
K/K)$-group. To construct our prolongation $G/X-S$ we 
shall not use the full
hypothesis of decency of $\scr G$---only that $\Gamma$ is 
finitely
presented\footnote{Ofer Gabber has sketched a construction 
of such a
``prolongation'' in the case where 
$G=\underline{\Aut}(V/K)$ for $V$ a
projective smooth variety over $K$ which requires only 
that $\Gamma$ be
finitely generated.}.  For simplicity of notation, let us 
first consider the
case where $\Gamma$  is constant rather than only locally 
constant, when viewed
as sheaf for the \'etale topology over $K$. Let $R\to 
F\to\Gamma\to 1$ be a
finite presentation of $\Gamma$, with $F$ a free group 
generated by elements
$x_1,\dotsc,x_n$ and $R$ a free group generated by 
elements $y_1,\dotsc,y_m$.
By some abuse of notation we may view the natural 
projection $G\to \scr G$ as
giving us a homomorphism $G\to\Gamma$, and let $G^\sharp$ 
denote the
fiber-product of $G\to\Gamma$ and $F\to\Gamma$, so that we 
have a commutative
diagram
$$
\CD
1 @>>> G^0 @>>> G^\sharp @>>> F @>>> 1\\
@.  @V=VV  @VVV  @VVV\\
1 @>>> G^0  @>>> G  @>>> \Gamma   @>>> 1
\endCD
$$
where the rows are exact. For any element $s$ in $F$ (or 
in $\Gamma)$ let
$G_s^0$ denote the $G^0$-bitorsor in $G^\sharp$ (or in 
$G)$ which is the full
inverse image of $s$. To say that a scheme $Y$ is a 
bitorsor for a group scheme
$H$ means that $Y$ has the structure of a left- and a 
right-\<$H$-torsor, these
actions commuting (cf. \cite{Grot3, Exp.~VII} and 
\cite{Br} for a treatment of
this notion).

First note that $G^0$ itself does prolong to a smooth 
group scheme over $X-S$
(for some finite $S)$. Fix such a prolongation; call it 
$G^0/X-S$. If $S'$ is
any finite set containing $S$, let $G^0/X-S'$ denote the 
restriction of
$G^0/X-S$ to $X-S'$. 

We shall now extend this prolongation of $G^0$ to a 
prolongation of $G^\sharp$
to a smooth locally algebraic group scheme over $X-S$ 
(possibly enlarging $S)$.
For this, consider the $G^0$-bitorsors $G_{x_1}^0,G_{x_2}^0,
\dotsc,G_{x_n}^0$ for the set
of generators $x_1,\dotsc,x_n$ of $F$, these all being 
bitorsors over the base
$\Spec K$. Each one of these bitorsors $G_{x_i}^0$ 
prolongs to a bitorsor for
$G^0/X-S_i$ for some finite set $S_i$ containing $S$. 
Replacing $S$ by the
union of the $S_i$'s, we have bitorsors $G_{x_i}^0/
X-S$ for $G^0/X-S$ (for
$i=1,\dotsc,n)$. For any word $w$ in $x_1,x_2,\dotsc,x_n$ 
and their inverses in
the free group $F$, letting $G_w^0/X-S$ be the appropriate
(``contraction-product''\footnote{This is the 
noncommutative version of the
Baer-sum: one contracts the product of bitorsors by the 
group action which is
the left-action on the right bitorsor and the right-action 
on the left
bitorsor.}) bitorsor for $G^0/X-S$ built from the 
$G_{x_i}^0/X-S$ one checks
that the disjoint union $\coprod G_w^0/X-S$ taken over all 
elements $w$ in $F$
(together with multiplication induced by the 
contraction-product construction)
constitutes a prolongation of $G^\sharp$ to a locally 
algebraic group scheme
over $X-S$ with the desired properties. Call it 
$G^\sharp/X-S$.

To ``descend'' this prolongation to a prolongation of $G$ 
one sees that it
suffices to choose trivializations of the $G^0/X-S$ 
bitorsors $G_{y_j}^0/X-S$
for each of the generators $y_j$ of the group of relations 
$R$, such that these
chosen trivializations are ``consistent'' in the sense 
that the induced
trivialization on $G_y^0/X-S$ for any $y\in R$ is well 
defined, i.e.,
independent of how $y$ is expressed in terms of the 
generators $y_j$. But since
$G^0$ is separated, ``consistency'' of these 
trivializations over $K$ implies
``consistency'' over $X-S$. We can find such a 
``consistent'' choice of
trivializations (for our generating set $y_j)$ over $K$, 
and then, after a
possible further enlargement of $S$ to a finite set of 
places of $K$, we may
extend this choice to $X-S$. 

We have then taken care of the case of constant $\scr G$, 
i.e., discrete
$\Gal(\overline K/K)$-groups $\Gamma$ with trivial Galois 
action. In general
though, we have that $\Gal(\overline K/K)$ acts on 
$\Gamma$ continuously, and
therefore, recalling the notation of \S16, the action 
factors through a finite
quotient, $\Delta=\Gal(L/K)$. We then begin the procedure 
as above, i.e., by
choosing a prolongation $G^0/X-S$ as before, then making 
the base change to
$Y-T$ where $Y$ is the spectrum of the ring of integers in 
$L$ and $T$ is the
inverse image of $S$, and continuing the construction over 
$Y-T$ appropriately
``equivariantly'' for the action of $\Gal(L/K)$, and then 
finally descending
back to $X-S$.\qed\enddemo

Let $G$ denote the extended group scheme given by the 
preceding lemma.  We know
by Theorem 4 that $\cyr{X}_S(G)$ is finite. Let us define 
an ``intermediate
group'' (call it $\cyr{X}(S;G))$ which will allow us to 
compare $\cyr{X}_S(G)$
and $\cyr{X}(G)$, namely,
$$
\cyr{X}(S,G)\coloneq \ker\left\{H^1(X-S,G)\to \prod_v 
H^1(K_v,G)\right\},
$$
where the summation over $v$ is over all places $v$ of 
$K$. Now, by definition,
$\cyr{X}(S,G)$ is contained in $\cyr{X}_S(G)$ and is 
consequently also finite.
Moreover, if $S\subset T$ is an inclusion of finite sets 
of places of $v$, the
natural mapping $\cyr{X}(S,G)\to \cyr{X}(T,G)$ is easily 
seen to be surjective.
Since $\cyr{X}(G)=\varinjlim\cyr{X}(T,G)$, the limit taken 
over the directed
system of all finite sets of places containing $S$, the 
corollary follows.\qed

Now let $V/K$ be a smooth projective variety. The 
automorphism group
$G=\underline{\Aut}(V/K)$ is a locally algebraic group 
over $K$. Let $\scr
G=G/G^0$ be its (discrete) group of connected components 
and $\Gamma(V)\coloneq
\scr G(\overline K)$ the group of $\overline K$-valued 
points of $\scr G$,
i.e., the group of connected components of the locally 
algebraic group
$G/\overline K=\underline{\Aut}(V/\overline K)$.

Very little seems to be known, in general, about the 
groups $\Gamma(V)$. Are
they finitely generated for all $V$? Finitely presented? 
decent? Are they all
arithmetic groups?

The group $\Gamma(V)$ admits a natural mapping (with 
finite kernel) to an
arithmetic group.\footnote{Under the hypothesis that 
$\Gamma(V)$ is finitely
generated, Ofer Gabber has sketched a proof that 
$\Gamma(V)$ can be embedded in
an arithmetic group.} To see this, one uses the natural 
representation of $G$
on $\NS(V)$, the N\'eron-Severi group of $V$ which we view 
as finitely
generated abelian locally constant group over $K$, whose 
full automorphism
group is an arithmetic group. One then proves that the 
subgroup $G_1$ of $G$
which is the identity on N\'eron-Severi is an algebraic 
group and, in
particular, has only a finite number of components.

Let us recall, briefly, how $G_1$ is analyzed. First, the 
action of $G_1$ on
$V$ gives us an action on $\Pic^0(V)$. If $L$ is any line 
bundle and $g\in
G_1$, then $[g_*L\otimes L^{-1}]=[(g^{-1})^*L\otimes 
L^{-1}]$ determines a class
in $\Pic^0(V)$, giving us a 1-cocycle $\varphi_L$ on the 
group scheme $G_1$
with values  in $\Pic^0(V)$ defined by $g\mapsto 
[g_*L\otimes L^{-1}]$. Fixing
$L$ to be the very ample line bundle on $V$ coming from a 
chosen projective
embedding, denote by $G_{1,L}\subset G_1$ the ``kernel'' 
of the 1-cocycle $
\varphi_L$. Define $G_{2,L}$ to be the group (representing 
the functor) of
paired automorphisms of $(V,L)$, i.e., automorphisms of 
$V$ together with
compatible automorphism of $L$. We have a surjection of 
locally algebraic
groups $G_{2,L}\to G_{1,L}$ (in fact, $G_{2,L}$ is a 
central extension of
$G_{1,L}$ with kernel isomorphic to $\Bbb G_m)$. But 
$G_{2,L}$ is clearly a
closed subgroup of the general linear group of 
automorphisms of $H^0(V,L)$. It
follows that $G_1$ is of finite type, and therefore its 
image in
$\Aut(\Pic^0(V))$ is finite.

Collecting things, then, we can say that 
$G=\underline{\Aut}(V/K)$ is a
``successive extension'' of a subgroup of $\Aut(\NS(V))$ 
by a subgroup of
$\Aut(\Pic^0(V))$ by a subgroup scheme in $\Pic^0(V)$ by 
a linear group.

\rem{Remark} The ``abelian variety part'' of $G$ is 
trivial if $\Pic^0(V)$
vanishes.
\endrem

Compare \cite{Ra} where the connected component containing 
the identity in
$\underline{\Aut}(V)$ is discussed. 

By \cite{BS, Theorem 2.6}, or by the more general 
``descent theory'' (cf.
\cite{Grot}) the set $\scr S(V/K)$ may be identified with 
$\cyr{X}(G)$. Hence,
we have

\thm{Corollary 2} Let $V/K$ be a projective variety such 
that the locally
constant group $\scr G(V/K)$ of components of 
$\underline{\Aut}(V/K)$ is
decent. Suppose, further, that the Tate-Shafarevich 
Conjecture holds \RM(for
all twists of abelian subvarieties of $\Pic^0(V/K)$ 
defined over $K)$. Then $
\scr S(V/K)$ is finite, i.e., the local-to-global 
principle holds for $V/K$, up
to finite obstruction. 
\ethm

\thm{Corollary 3} Let $V/K$ be a projective variety such 
that $\Pic(V/\overline
K)=\Bbb Z$. Then $\scr S(V/K)$ is finite, i.e., the 
local-to-global principle
holds for $V/K$, up to finite obstruction.
\ethm

\demo{Proof} This follows directly from Corollary 2. For 
if $\Pic(V/\overline
K)=\Bbb Z$, then $\scr G(V/K)$ is a finite locally 
constant group scheme and,
therefore, is decent. Also, since $\Pic^0(V/K)$ vanishes, 
the hypothesis
involving finiteness of the Tate-Shafarevich group of 
twists of subabelian
varieties of $\Pic^0(V)$ is certainly satisfied.\qed\enddemo

\thm{Corollary 4} Let $V/K$ be a smooth variety of 
dimension $\ge 3$ which is a
hypersurface \RM(or more generally, a complete 
intersection\RM) in projective
space. Then $\scr S(V/K)$ is finite, i.e., the 
local-to-global principle holds
for $V/K$, up to finite obstruction.
\ethm

\demo{Proof} For such varieties, $\Pic(V/\overline K)=\Bbb 
Z$ \cite{Grot2,
Corollary 3.7 of Exp.~XII} and therefore Corollary 3 
applies.\qed\enddemo

I am thankful to O. Gabber, J. Harris, and J.-P. Serre who 
told me the proof of
the following proposition, which is a result due to Jordan 
\cite{J} and which
implies that, for smooth hypersurfaces $V$ of dimension 
$\ge 3$ and of degree
$\ge 3$, $\underline{\Aut}(V)$ is, in fact, finite. 
Therefore, for such
varieties the conclusion of Corollary 4 follows directly.

\thm{Proposition} Let $n\ge 2$, and let $V$ be a smooth 
hypersurface in 
$\bold P^n/
\Bbb C$ of degree $\ge 3$. Let $\Phi$ denote the subgroup 
of $\PGL_{n+1}(\Bbb
C)$ which stabilizes $V$. Then $\Phi$ is finite.
\ethm

\demo{Proof} Let $V$ be defined by the homogeneous form 
$F(X_0,\dotsc,X_n)$ of
degree $d$. The group $\Phi$ is algebraic; it suffices to 
show that its Lie
algebra vanishes. Suppose, then, that $A=(a_{ij})$ 
$(i,j=0,\dotsc,n)$ is an
$(n+1)\times(n+1)$ matrix in the Lie algebra of $\GL_{n+
1}$ which stabilizes
$V$, i.e., such that
$$
\sum_{i,j}a_{ij}\bfcdot X_j\bfcdot\partial F/\partial 
X_i=\lambda\bfcdot F\quad
\text{for }\lambda\in\Bbb C^*\.
$$

Writing $\lambda\bfcdot F=\lambda\bfcdot d^{-1}\sum 
\delta_{ij}\bfcdot
X_j\bfcdot\partial F/\partial X_i$ where $(\delta_{ij})$ 
is the identity
matrix, we have
$$
\sum_{i,j}(a_{ij}-\lambda d^{-1}\delta_{ij})X_j\bfcdot 
\partial F/\partial
X_i=0\.\tag 2
$$

Since $V$ is smooth, the forms $\partial F/\partial X_i$ 
only have the origin
in $\Bbb C^{n+1}$ as a common intersection, i.e., the 
sequence 
$$
(\partial F/\partial X_0, \partial F/\partial 
X_1,\dotsc,\partial F/\partial
X_n)
$$
is a ``regular sequence''; therefore, the ideal of 
relations of this (regular)
sequence is generated by the evident ones 
$$
\partial F/\partial X_i\bfcdot \partial F/\partial 
X_j-\partial F/\partial
X_j\bfcdot \partial F/\partial X_i=0\.
$$

In other words, the kernel $\scr R$ of the homomorphism of 
the $\Bbb
C[X_0,\dotsc,X_n]$-algebra $\Bbb C[X_0,\dotsc, 
X_n,Y_0,\dotsc,Y_n]$ to
$\Bbb C[X_0,\dotsc,X_n]$ which sends the variable $Y_j$ to 
$\partial F/\partial
X_j$ $(j=0,\dotsc,n)$ is generated by the elements
$$
\partial F/\partial X_j\bfcdot Y_i-\partial F/\partial 
X_i\bfcdot Y_j
\qquad (i,j=0,\dotsc,n)\.
$$
Since these elements are homogeneous of degree $d-1\ge 2$ 
in the $X_k$'s, it
follows that if $\rho_j(X_0,\dotsc,X_n)$ are linear forms 
in the $X_k$'s for
$j=0,\dotsc,n$, such that 
$$
\sum_j\rho_j(X_0,\dotsc,X_n)\bfcdot\partial F/\partial X_j=0
$$
then each of the forms $\rho_j(X_0,\dotsc,X_n)$ vanish 
identically. Applying
this to (2) we have that $a_{ij}=\lambda 
d^{-1}\delta_{ij}$, i.e., that $A$
gives $0$ in the Lie algebra of $\PGL_{n+1}$.\qed\enddemo

\rem{Remarks} Applying this proposition to smooth 
hypersurfaces of dimension
$\ge 3$ and degree $\ge 3$, and using (as in the proof of 
Corollary 4) the fact
that, for such varieties, $\Pic(V)=\Bbb Z$, we see that 
$\underline{\Aut}(V)$
is finite. See the discussion of this proposition on 
\cite{B, pp.~41--42}.
\endrem

\ex{Exercise (Serre)} What happens when the base field is 
algebraically closed
of characteristic $p$, and $p$ divides the degree $d$ of 
$V$?
\medskip
What analogous results are there for smooth hypersurfaces 
in multiprojective
spaces?
\endex

\subheading{18. The mapping $\cyr{X}(A/K)\to\scr S(A/K)$} 
If $A$ is an abelian
scheme over a base $S$, to give an automorphism of $A/S$ 
as abelian scheme is
equivalent to giving an automorphism of $A$ viewed merely 
as scheme over $S$,
which fixes the zero-section. Let $\Aut_\bullet(A)$ denote 
the group of such
automorphisms. Fixing $A/K$, an abelian variety, let
$\underline{\Aut}_\bullet(A/K)$ denote the locally 
algebraic group of
automorphisms of $A$ preserving zero-section, i.e., of the 
abelian variety $A$.
Let $\scr A$ denote the set of isomorphism classes of 
abelian varieties $A'/K$
which are isomorphic to $A/K_v$ as abelian varieties over 
$K_v$, for all places
$v$ of $K$. The set $\scr A$ is finite \cite{BS, Corollary 
7.11}. In fact, $
\scr A$ can be identified with the set
$\cyr{X}(\underline{\Aut}_\bullet(A/K))$, and, since
$\underline{\Aut}_\bullet(A/K)$ is a locally constant 
group of ``type ALA'' in
the terminology of \cite{BS}, the required finiteness 
follows.

\thm{Theorem 5} The set $\scr S(A/K)$ has a natural 
partition into a finite
number of \RM(disjoint\RM) subsets indexed by $\scr A$,
$$
\scr S(A/K)=\coprod_{\alpha\in\scr A}\scr S(A/K)_\alpha,
$$
and, for any $\alpha\in \scr A$, if $A'/K$ is an abelian 
variety in the
isomorphism class $\alpha$, then there is a natural {\it 
surjection\/}
$$
\cyr{X}(A'/K)\to \scr S(A/K)_\alpha\tag 3
$$
which identifies $\scr S(A/K)_\alpha$ with the orbit-space 
of $\cyr{X}(A'/K)$
under the action of the group of automorphisms of $A'$ 
defined over $K$.
\ethm

\demo{Proof} First, let $V$ be a ``companion'' to $A$, 
i.e., a projective
variety over $K$, isomorphic to $A$ over the completions 
$K_v$ where $v$ ranges
over all places of $K$. For a model of the ``Albanese 
variety'', $\A_V$, of $V$
let us take $\A_V=\Pic^0(\Pic^0(V))$. Thus $\A_V$ is an 
abelian variety over
$K$ isomorphic, as abelian variety, to $A$ over every 
completion of $K$ and
therefore the isomorphism class of the abelian variety 
$\A_V$ over $K$ is an
element of $\scr A$. There is a natural principal 
homogeneous action of $\A_V$
on $V$ (defined over $K)$. Since $V$ is also a companion 
to $\A_V$ over $K$, it
(taken with its structure as principal homogeneous space 
for $\A_V)$ determines
an element in $\cyr{X}(\A_V)$. 

Next, for $\alpha\in\scr A$, choose an abelian variety 
$A'/K$ in the
isomorphism class $\alpha$, and define $\scr 
S(A/K)_\alpha$ to be the subset of
$\scr S(A/K)$ represented by varieties $V/K$ such that 
$\A_V\cong A'$ (the
isomorphism being as abelian varieties over $K)$. This 
clearly gives a
partition of $\scr S(A/K)$. For each $\alpha\in\scr A$, we 
have the required
surjection (3), and this mapping visibly factors through 
the quotient to the
orbit-space under the group of $K$-automorphisms of $A'$.

To conclude the theorem, we consider two elements in 
$\cyr{X}(A'/K)$ which
represent the same element in $\scr S(A/K)_\alpha$. 
Equivalently, on the same
projective variety $V$, we are given two principal 
homogeneous space structures
under the abelian variety $A'$. We must show that one can 
bring one of these
structures to the other by a $K$-automorphism $\alpha$ of 
$A'$. But this is
easy: We get a natural morphism over $K$ from the 
projective variety $V$ to the
(locally constant) group $\underline{\Aut}_\bullet(A'/K)$ 
by the rule which
associates to each point $v$ of $V$ the mapping 
$\varphi_v\:A'\to A'$ given as
follows: to any $a\in A'$, the translate of $v$ by $a$ via 
the first principal
homogeneous structure is equal to the translate of $v$ by 
$\varphi_v(a)$ via
the second principal homogeneous structure. Since $V$ is 
connected and
$\underline{\Aut}_\bullet(A'/K)$ is locally constant, the 
mapping $v\mapsto 
\varphi_v$ is, in fact, constant, and therefore its image 
is an automorphism $
\alpha$ of $A'$, defined over $K$.\qed\enddemo

When is the mapping $\cyr{X}(A'/K)\to\scr S(A/K)$ 
surjective? From the above
analysis, this will happen if and only if
$\cyr{X}(\underline{\Aut}_\bullet(A/K))$ is trivial. Here 
are two instances:

\thm{Corollary 1} $\cyr{X}(A'/K)\to\scr S(A/K)$ is 
surjective if either\RM:
\roster
\item "(a)" there is a finite cyclic extension $L/K$ such 
that every
automorphism of $A$ over $\overline K$ is already defined 
over $L$, or
\item "(b)" \<$A=E$ is an elliptic curve.
\endroster
\ethm

\demo{Proof} Since $\underline{\Aut}_\bullet(A/K)$ is 
discrete, case (a)
follows immediately from the corollary to Example (1) in 
\S15. Case (b) follows
from case (a), for if $A$ is an elliptic curve over a 
number field $K$, the
splitting field over $K$ for the action of $\Gal(\overline 
K/K)$ on
$\underline{\Aut}_\bullet(A)$ is of degree $\le 2$.\qed\enddemo

\thm{Corollary 2} For $K$ a fixed number field, the 
Tate-Shafarevich Conjecture
holds for all abelian varieties over $K$ if and only if 
$\scr S(A/K)$ is finite
for all abelian varieties $A$ over $K$.
\ethm

\demo{Proof} Recall that $\scr A$ is a finite set. In view 
of the above
theorem, then, it suffices to show that, for any abelian 
variety,
$A/K$, $\cyr{X}(A/K)$ is finite if and only if its image 
in $\scr S(A/K)$ is
finite, or, equivalently, if and only if there are only a 
finite number of
orbits of $\cyr{X}(A/K)$ under the action of the 
$K$-automorphism group of
$A$. But, as Tate remarked to me, this is easy to see: 
there are only a finite
number of elements in $\cyr{X}(A/K)$ of any fixed order, 
and automorphisms of
groups preserve the order of elements.\qed\enddemo

\subheading{19. The local-to-global principle for 
quadrics} I am thankful to
Colliot-Th\'el\`ene for showing me this proof of the 
local-to-global principle
for any smooth quadric $V$ over a number field $K$.

Let $W/K$ be a smooth projective variety which is {\it 
potentially quadric\/}
in the sense that its base change to $\overline K$ is 
isomorphic to a quadratic
hypersurface in projective space. Putting $\overline 
W=W\times_{\Spec
K}\Spec\overline K$ we have that $\Pic(\overline W)$ is 
isomorphic either to
$\Bbb Z$ or $\Bbb Z\oplus \Bbb Z$ (the latter only if 
$\overline W$ is of
dimension 2).  Call an element $\xi$ of $\Pic(\overline 
W)$ {\it of quadric
type\/} if a line bundle $\overline{\scr L}$ over 
$\overline W$ representing
$\xi$ has the property that $\overline{\scr L}$ is very 
ample and the embedding
of $\overline W$ in projective space associated to the 
complete linear system
of $\overline{\scr L}$ identifies $\overline W$ with a 
quadric hypersurface
over $\overline K$. There is only one element $\xi$ ``of 
quadric type'' in
$\Pic(\overline W)$, and that element $\xi$ is fixed under 
the action of
$\Gal(\overline K/K)$. We have the exact sequence
$$
0\to\Pic(W)\to\Pic(\overline W)^{\Gal(\overline 
K/K)}\underset{\delta_K}\to \to
\Br(K)\tag 4
$$
where $\Br(K)$ is the Brauer group of $K$ and $\delta_K$ 
is the natural
``coboundary'' mapping whose definition is as follows. For 
a line bundle
$\overline l$ on $\overline W$ whose isomorphism class is 
fixed by
$\Gal(\overline K/K)$ and for each 
$\sigma\in\Gal(\overline K/K)$ choose an
isomorphism $\iota_\sigma\:\overline l\cong\overline 
l^\sigma$. Denote by $c$
the 2-cocycle on $\Gal(\overline K/K)$ with coefficients 
in $\overline K^*$
defined by the formula 
$c(\sigma,\tau)=(\iota_{\sigma\tau})^{-1}\bfcdot(\iota_
\sigma)^\tau\bfcdot\iota_\tau$. Then if $[\overline l]$ 
denotes the isomorphism
class of $\overline l$ in $\Pic(\overline 
W)^{\Gal(\overline K/K)}$,
$\delta_K([\overline l])$ is defined to be the cohomology 
class of $c$ in the
Brauer group $\Br(K)=H^2(\Gal(\overline K/K);\overline 
K^*)$; this cohomology
class is seen to be dependent only upon the isomorphism 
class of $[l]$ and
independent of the choice of the $\iota_\sigma$'s. 

The element $\delta_K(\xi)\in\Br(K)$ may be thought of as 
the obstruction to
$W/K$ being isomorphic to a quadric hypersurface over $K$, 
for its vanishing is
a necessary and sufficient condition for the line bundle 
$\overline{\scr L}$
corresponding to $\xi$ on $\overline W$ to come from a 
line bundle $\scr L$ on
$W/K$ whose linear system would identify $W$ with a 
quadric hypersurface in
projective space over $K$.

If $V$ is a smooth quadric over $K$ and $W$ a companion to 
$V$ over $K$, then
$V$ is isomorphic to $W$ over $\overline K$, so $W$ is 
potentially quadric.
Since $V/K_v$ is isomorphic $W/K_v$ for all places $v$ of 
$K$, $W$ is, in fact,
quadric over $K_v$ for all $v$; therefore, it follows from 
the version of (4)
over $K_v$ that $\delta_{K_v}(\xi)=0$ in $\Br(K_v)$ for 
all $v$. Since the
global Brauer group injects into the direct sum of the 
local Brauer groups,
$\delta_K(\xi)=0$; therefore, the companion $W$ is a 
quadric variety over
$K$.\footnote{A similar argument shows that a smooth 
projective variety $W/K$
of dimension $\ge 3$ which is a companion to a 
hypersurface of degree $d$
defined over $K$ is itself isomorphic over $K$ to a 
hypersurface of degree $d$.}
Now let $q_V$ and $q_W$ denote the two quadratic forms 
over $K$ corresponding
to $V$ and $W$. Since $V$ and $W$ are companions, it 
follows (from the argument
alluded to in the footnote to \S2 of Part I) that the 
corresponding forms
$q_{V,v}$ and $q_{W,v}$ over the completions $K_v$ for 
each place $v$ of $K$
are similar. Using \cite{O} we then see that $q_V$ and 
$q_W$ are similar over
$K$.\qed

\heading Acknowledgment\endheading
I am thankful to M. Artin, P. Diaconis, J.-L. 
Colliot-Th\'el\`ene, O. Gabber,
J.~Harris, N.~Katz, Y.~Nisnevich, K.~Rubin, and J.-P.~Serre
for their generous
help and to the audiences at the MSRI\footnote{At the 
conference in June 1992
commemorating the 10th anniversary of the MSRI and in 
honor of Irving Kaplansky
on the occasion of his retirement as director.}, the 
University of Toronto, and
Penn State\footnote{The Weisfeiler Lecture, November 
1992.}, for their
comments, suggestions, and patience during my talks on 
material corresponding
to Part I of this article. 

\Biblio
\widestnumber\key{C-T--S-D}
\ref
\key A1
\by M. Artin
\paper Algebraic spaces
\inbook Yale Math. Monographs
\vol 3
\publ Yale Univ. Press
\publaddr New Haven, CT
\yr 1976
\endref

\ref
\key A2
\by M. Artin
\book N\'eron models
\bookinfo Arithmetic Geometry
\eds G. Cornell and J. Silverman
\publ Springer-Verlag
\publaddr New York
\yr 1986
\pages 213--230
\endref

\ref
\key BD
\by M. Bertolini and H. Darmon
\paper Kolyvagin\RM's descent and Mordell-Weil groups over 
ring class fields
\jour J. Reine Angew. Math.
\vol 412
\yr 1990
\pages 63--74
\endref

\ref
\key B
\by A. Borel 
\book Introduction aux groupes arithm\'etiques
\publ Hermann
\publaddr Paris
\yr 1969
\endref

\ref
\key BS
\by A. Borel and J.-P. Serre 
\paper Th\'eor\`emes de finitude en cohomologie galoisienne
\jour Comment. Math. Helv. 
\vol 39 
\yr 1964
\pages 111--164
\finalinfo
{\it Oeuvres. Collected papers\/}, vol. II 1959-1968, 
Springer-Verlag,
New York, 1983, pp. 362--415
\endref

\ref
\key BLR
\by S. Bosch, W. L\"utkebohmert, and M. Raynaud
\paper N\'eron models
\inbook Ergeb. Math. Grenzgeb., bd. 21
\publ Springer-Verlag
\publaddr New York
\yr 1990
\endref

\ref
\key Br
\by L. Breen
\paper Bitorseurs et cohomologie non ab\'elienne
\inbook The Grothendieck Festschrift,
{\rm vol. II}
(A collection of articles written in honor of the 60th 
birthday of
Alexander Grothendieck) 
\publ Birkh\"auser
\publaddr Basel
\yr 1991
\pages 399--476
\endref

\ref
\key Ca1
\by J. W. S. Cassels 
\book Rational quadratic forms
\publ Academic Press
\publaddr New York
\yr 1978
\endref

\ref
\key Ca2
\bysame
\paper Lectures on elliptic curves
\inbook London Math. Soc. Stud. Texts
\vol 24
\publ Cambridge Univ. Press
\publaddr London and New York
\yr 1991
\endref

\ref
\key Ch
\by V. Chernousov
\paper On the Hasse principle for groups of type $E_8$
\jour Soviet Math. Dokl. 
\vol 39
\yr 1989 
\pages 592--596
\endref

\ref
\key Coa
\by J. Coates
\paper Elliptic curves with complex multiplication and 
Iwasawa theory
\jour Bull. London Math. Soc. 
\vol 23
\yr 1991
\pages 321--350
\endref

\ref
\key Coh
\by P. J. Cohen
\paper Decision procedures for real and $p$-adic fields
\jour Comm. Pure Appl. Math.
\vol 22
\yr 1969
\pages 131--151
\endref

\ref
\key C-T
\by J.-L. Colliot-Th\'el\`ene 
\paper Les grands th\`emes de Fran\c cois Ch\^atelet
\jour Enseign. Math. 
\vol 34
\yr 1988
\pages 387--405
\endref

\ref
\key C-T--S-D
\by J.-L. Colliot-Th\'el\`ene and P. Swinnerton-Dyer
\book Hasse principle and weak approximation for pencils 
of Severi-Brauer and
similar varieties
\finalinfo preprint, 1992
\endref

\ref
\key D1
\by H. Darmon
\book Euler systems and refined conjectures of Birch 
Swinnerton-Dyer type
\bookinfo Proc. Conf. $p$-adic Birch Swinnerton-Dyer 
Conjectures 
(Boston Univ., August
1991)
\toappear
\endref

\ref
\key D2
\bysame
\paper A refined conjecture of Mazur-Tate type
\jour Invent. Math.
\vol 110
\yr 1992
\pages 123--146
\endref

\ref
\key Du
\by W. Duke
\paper Hyperbolic distribution problems and half-integral 
weight Maass forms
\jour Invent. Math.
\vol 92
\yr 1988
\pages 73--90
\endref

\ref
\key D--S-P
\by W. Duke and R. Schulze-Pillot
\paper Representation of integers by positive ternary 
quadratic forms and
equidistribution of lattice points on ellipsoids
\jour Invent. Math.
\vol 99
\yr 1990
\pages 49--57
\endref

\ref
\key F
\by M. Flach
\paper A finiteness theorem for the symmetric square of an 
elliptic curve 
\jour Invent. Math. 
\vol 109 
\yr 1992
\pages 307--327
\endref

\ref
\key Gi
\by J. Giraud
\paper Cohomologie non ab\'elienne
\inbook Grundlehren Math. Wissen., bd. 179
\publ Springer-Verlag
\publaddr New York
\yr 1971
\endref

\ref
\key GF
\by E. Golubeva and O. Fomenko
\paper Application of spherical functions to a certain 
problem in the theory of
quadratic forms
\jour J. Soviet Math.
\vol 38
\yr 1987
\pages 2054--2060
\endref

\ref
\key Gros
\by B. Gross
\paper Kolyvagin\RM's work on modular elliptic curves
\inbook $L$-functions and Arithmetic 
\eds J. Coates and M. J. Taylor
\bookinfo Proc. Durham 1989, London Math. Soc. Lecture 
Note  Ser.
\vol 153
\publ Cambridge Univ. Press
\publaddr London and New York
\yr 1991
\endref

\ref
\key Grot1
\by A. Grothendieck
\paper Techniques de descente et th\'eor\`emes 
d\,\RM'existence en 
g\'eom\'etrie
alg\'ebrique
\inbook S\'em. Bourbaki
\vol 12 (59/60)
\publ Benjamin Inc.
\publaddr New York and Amsterdam
\yr 1966
\endref

\ref
\key Grot2
\by A. Grothendieck (ed.)
\paper Cohomologie locale des faisceaux coh\'erents et 
Th\'eor\`emes de
Lefschetz locaux et globaux {\rm (SGA 2)}
\inbook Advanced Studies in Pure Mathematics 
\vol 2
\publ North-Holland
\publaddr Amsterdam 
\yr 1968
\endref

\ref
\key Grot3
\by A. Grothendieck (ed.)
\paper Groupes de monodromie en g\'eom\'etrie alg\'ebrique
{\rm (SGA 7 I)}
\inbook S\'eminaire de G\'eom\'etrie Alg\'ebrique du 
Bois-Marie 1967--1969
\bookinfo Lecture Notes in Math.
\vol 288
\publ Springer-Verlag
\publaddr New York
\yr 1972
\endref

\ref
\key Hard
\by G. Harder
\paper \"Uber die Galoiskohomologie halbeinfacher 
Matrizengruppen
\jour Math. Z. 
\vol 90
\yr 1965
\pages 404--428
\moreref
\vol 92
\yr 1966
\pages 396--415
\endref

\ref
\key Harr
\by J. Harris
\paper Algebraic geometry \RM(a first course\RM)
\inbook  Graduate Texts in Math.
\vol 133
\publ Springer-Verlag
\publaddr New York
\yr 1992
\endref

\ref
\key Hart
\by R. Hartshorne 
\book Algebraic geometry
\publ Springer-Verlag
\publaddr New York
\yr 1977
\endref

\ref
\key Has
\by H. Hasse
\paper \"Uber die \"Aquivalenz quadratischer Formen im 
K\"orper der rationalen
Zahlen
\jour J. Reine Angew. Math. 
\vol 152
\yr 1923
\pages 205--224
\endref

\ref
\key I
\by H. Iwaniec
\paper Fourier coefficients of modular forms of half 
integral weight
\jour Invent. Math. 
\vol 87
\yr 1987
\pages 385--401
\endref

\ref
\key J
\by C. Jordan
\paper M\'emoire sur l\RM'equivalence des formes
\jour J. \'Ecole Polytech 
\vol XLVIII
\yr 1880
\pages 112--150
\endref

\ref
\key KO
\by S. Kobayashi and T. Ochiai
\paper Mappings into compact complex manifolds with 
negative first Chern class
\jour J. Math. Soc. Japan
\vol 23
\yr 1971
\pages 136--148
\endref

\ref
\key K1
\by V. A. Kolyvagin
\paper Finiteness of $E(\Bbb Q)$ and $\cyr{X}(\Bbb Q)$ for 
a subclass of Weil
curves
\jour Izv. Akad. Nauk SSSR Ser. Mat.
\vol 52
\yr 1988
\pages 522--540
\transl English transl.\jour Math USSR-Izv. 
\vol 32
\yr 1989
\pages 523--542
\endref

\ref
\key K2
\bysame
\paper On the Mordell-Weil group and the Shafarevich-Tate 
group of Weil
elliptic curves
\jour Izv. Akad. Nauk SSSR Ser. Mat.
\vol 52
\yr 1988
\pages 1154--1180
\endref

\ref
\key K3
\bysame 
\paper Euler systems
\inbook The Grothendieck Festschrift
\bookinfo vol. II (A collection of articles written in 
honor of the 60th
birthday of Alexander Grothendieck)
\publ Birkh\"auser
\publaddr Basel
\yr 1991
\pages 435--483
\endref

\ref
\key KL
\by V. A. Kolyvagin and D. Logachev
\paper Finiteness of the Shafarevich-Tate group and group 
of rational points
for some modular abelian varieties
\jour Algebra i Analiz
\vol 1
\yr 1989
\pages 171--196
\transl English transl.,\nofrills
\jour Leningrad Math. J.
\vol 1
\yr 1990
\pages 1229--1253
\endref

\ref
\key Kne
\by M. Kneser
\paper Lectures on Galois cohomology of classical groups
\inbook Tata Inst. Fund. Res. Lectures on Math. and Phys.
\publ Tata Inst. Fund. Res.
\publaddr Bombay
\yr 1969
\endref

\ref
\key Knu
\by D. Knutsen 
\paper Algebraic spaces
\inbook Lecture Notes in Math.
\vol 203
\publ Springer-Verlag
\publaddr New York
\yr 1971
\endref

\ref
\key Ma1
\by B. Mazur 
\paper Rational points of abelian varieties with values in 
towers of number
fields
\jour Invent. Math. 
\vol 18
\yr 1972
\pages 183--266
\endref

\ref
\key Ma2
\bysame
\paper Number theory as gadfly
\jour Amer. Math. Monthly
\vol 98
\yr 1991
\pages 593--610
\endref

\ref
\key Mat
\by Y. Matijasevic
\paper Enumerable sets are diophantine
\jour  Dokl. Akad. Nauk SSSR
\vol 191
\yr 1970
\transl English transl. 
\jour Soviet Math. Dokl. 
\vol 11
\yr 1970
\pages 354--358
\endref

\ref
\key Mi1
\by J. Milne 
\book The conjectures of Birch and Swinnerton-Dyer for 
constant abelian
varieties over function fields
\bookinfo Thesis, Harvard Univ.
\yr 1967
\endref

\ref
\key Mi2
\bysame \paper The Tate-Shafarevich group of a constant 
abelian variety
\jour Invent. Math. 
\vol 5 
\yr 1968
\pages 63--84
\endref

\ref
\key Mi3
\bysame
\paper On a conjecture of Artin and Tate
\jour Ann. of Math. (2)
\vol 102
\yr 1975
\pages 517--533
\endref

\ref
\key Mi4
\bysame
\book Arithmetic duality theorems
\publ Academic Press
\publaddr New York
\yr 1986
\endref

\ref
\key Mo
\by L. J. Mordell
\book Diophantine equations
\publ Academic Press
\publaddr New York
\yr 1969
\endref

\ref
\key Ne1
\by J. Nekov\'a\v r
\paper Kolyvagin\RM's method for Chow groups of Kuga-Sato 
varieties
\jour Invent Math.
\vol 107
\yr 1992
\pages 99--125
\endref

\ref
\key Ne2
\bysame 
\book $p$-Adic methods in arithmetic
\bookinfo Center for Pure and Appl. Math.
\publ Univ. of California
\publaddr Berkeley, preprint
\yr 1992
\endref

\ref
\key N\'e
\by A. N\'eron
\paper Mod\`eles minimaux des vari\'et\'es ab\'eliennes 
sur les corps locaux et
globaux
\jour Inst. Hautes \'Etudes Sci. Publ. Math.
\vol 21
\yr 1964
\endref

\ref
\key Nik1
\by V. Nikulin
\paper On factor groups of the automorphism groups of 
hyperbolic forms modulo
subgroups generated by $2$-reflections
\jour Soviet Math. Dokl.
\vol 20
\yr 1979
\pages 1156--1158
\endref

\ref
\key Nik2
\bysame
\paper Quotient-groups of the automorphisms of hyperbolic 
forms by subgroups
generated by $2$-reflections, algebro-geometric applications
\inbook Current Problems in Math., vol. 18
\bookinfo Akad. Nauk SSSR, Usesoyuz. Inst. Nauchn. Tekhn. 
Informatsii, Moscow
\yr 1981
\pages 3--114 
\transl English transl. 
\jour J. Soviet Math. 
\vol 22
\yr 1983
\issue 4
\endref

\ref
\key Nis1
\by Y. Nisnevich
\book Etale cohomology and arithmetic of semi-simple groups
\bookinfo Ph.D. thesis
\publ Harvard Univ.
\yr 1982
\endref

\ref
\key Nis2
\bysame 
\book On certain arithmetic and cohomological invariants 
of semi-simple
groups
\publ Johns Hopkins Univ., preprint
\yr 1989
\endref

\ref
\key O
\by T. Ono
\paper Arithmetic of orthogonal groups
\jour J. Math. Soc. Japan
\vol 7
\yr 1955
\pages 79--91
\endref

\ref
\key P
\by B. Perrin-Riou
\paper Travaux de Kolyvagin et Rubin
\jour S\'em. Bourbaki, exp. 717
\yr 1989/90
\moreref
\jour Ast\'erique 
\vol 189--190
\yr 1990
\pages 69--106
\endref

\ref
\key P-S--S
\by I. Piatetski-Shapiro and I. Shafarevich
\paper A Torelli theorem for algebraic surfaces of type $K3$
\jour Math. USSR-Izv. 
\vol 5
\yr 1971
\pages 547--588
\endref

\ref
\key Pl
\by V. P. Platonov 
\paper The arithmetic theory of algebraic groups
\jour Uspekhi Mat. Nauk
\vol 37
\yr 1982
\pages 3--54
\transl English transl., \nofrills
\jour Russian Math Surveys
\vol 37
\yr 1982
\pages 1--62
\endref

\ref
\key Ra
\by C. P. Ramanujam
\paper A note on automorphism groups of algebraic varieties
\jour Math. Ann. 
\vol 19
\yr 1964
\pages 25--33
\endref

\ref
\key R1
\by K. Rubin
\paper Tate-Shafarevich groups and $L$-functions of 
elliptic curves with
complex multiplication
\jour Invent. Math.
\vol 89
\yr 1987
\pages 527--560
\endref

\ref
\key R2
\bysame
\paper On the main conjecture of Iwasawa theory for 
imaginary quadratic fields
\jour Invent. Math.
\vol 93
\yr 1988
\pages 701--713
\endref

\ref
\key R3
\bysame
\paper Kolyvagin\RM's system of Gauss sums
\inbook Arithmetic Algebraic Geometry 
\bookinfo Progr. Math.
\vol 89
\eds van der Geer, Oort, and Steenbrink
\publ Birkh\"auser
\publaddr Boston, MA
\yr 1991
\pages 309--324
\endref

\ref
\key R4
\bysame
\paper The main conjecture 
\paperinfo Appendix to S. Lang's Cyclotomic Fields I and 
II, combined second
ed.
\inbook Graduate Texts in Math.
\vol 121
\publ Springer-Verlag
\publaddr New York
\pages 397--429
\yr 1990
\endref

\ref
\key R5
\bysame 
\paper The {\rm ``}main conjectures{\rm ''} of Iwasawa 
theory for imaginary
quadratic fields 
\jour Invent. Math.
\vol 103
\yr 1991
\pages 25--68
\endref

\ref
\key R6
\bysame 
\paper The one-variable main conjecture for elliptic 
curves with complex
multiplication
\inbook $L$-functions and Arithmetic
\bookinfo Proc. Durham 1989, London Math. Soc. Lecture 
Note Ser., vol. 153
\eds J. Coates and M. J. Taylor
\publ Cambridge Univ. Press
\publaddr London and New York
\yr 1991
\endref

\ref
\key R7
\bysame
\paper Stark units and Kolyvagin\RM's {\rm ``}Euler 
systems{\rm ''}
\jour J. Reine Angew. Math.
\vol 425
\yr 1992
\pages 141--154
\endref

\ref
\key R8
\bysame 
\paper The work of Kolyvagin on the arithmetic of elliptic 
curves
\inbook Arithmetic of Complex Manifolds
\bookinfo Lecture Notes in Math.
\vol 1399
\publ Springer-Verlag
\publaddr New York
\yr 1989
\pages 128--136
\endref

\ref
\key Se1
\by J.-P. Serre
\paper G\'eom\'etrie alg\'ebrique et g\'eom\'etrie 
analytique
\jour Ann. Inst. Fourier (Grenoble)
\vol 6
\yr 1956
\pages 1--42
\endref

\ref
\key Se2
\bysame 
\paper Cohomologie galoisienne
{\rm (3rd. ed.)}
\inbook Lecture Notes in Math.
\vol 5
\publ Springer-Verlag
\publaddr New York
\yr 1965
\endref

\ref
\key Se3
\bysame 
\book Local fields
\publ Springer-Verlag
\publaddr New York
\yr 1979
\endref

\ref
\key Si
\by C. Siegel 
\paper Equivalence of quadratic forms
\jour Amer. J. Math.
\vol 63
\yr 1941
\pages 658--680
\moreref
\book Carl Ludwig Siegel Gesammelte Abhandlungen, {\rm bd. 
II}
\publ Springer-Verlag
\publaddr New York
\yr 1966
\pages 217--239
\endref

\ref
\key St
\by N. Stephens
\paper The diophantine equation $X^3+Y^3=DZ^3$ and the 
conjectures of Birch and
Swinnerton-Dyer
\jour J. Reine Angew. Math.
\vol 231
\yr 1968
\pages 121--162
\endref

\ref
\key Ta
\by J. Tate
\paper The arithmetic of elliptic curves
\jour Invent. Math.
\vol 23
\yr 1974
\pages 179--206
\endref

\ref
\key Th
\by F. Thaine
\paper On the ideal class groups of real abelian number 
fields
\jour Ann. of Math. (2) 
\vol 128
\yr 1988
\pages 1--18
\endref

\ref
\key W
\by A. Weil
\book Sur l\RM'analogie entre les corps de nombres 
alg\'ebriques et les corps
de fonctions alg\'ebriques {\rm (1939)}
\bookinfo Oeuvres Sci. I
\publ Springer-Verlag
\publaddr New York
\pages 236--240
\yr 1980
\endref
\endRefs

\enddocument